\documentclass{article}%
\usepackage{amsfonts}
\usepackage{graphicx}
\usepackage{amsmath}
\usepackage{amssymb}%
\providecommand{\U}[1]{\protect\rule{.1in}{.1in}}
\newtheorem{theorem}{Theorem}

\newtheorem{corollary}[theorem]{Corollary}

\newtheorem{lemma}[theorem]{Lemma}

\begin{document}

\title{A Spectral Method for the Eigenvalue Problem for Elliptic Equations}
\author{Kendall Atkinson\\Departments of Mathematics \& Computer Science \\The University of Iowa
\and Olaf Hansen\\Department of Mathematics \\California State University San Marcos}
\maketitle

\begin{abstract}
Let $\Omega$ be an open, simply connected, and bounded region in
$\mathbb{R}^{d}$, $d\geq2$, and assume its boundary $\partial\Omega$ is
smooth. Consider solving the eigenvalue problem $Lu=\lambda u$ for an elliptic
partial differential operator $L$ over $\Omega$ with zero values for either
Dirichlet or Neumann boundary conditions. We propose, analyze, and illustrate
a `spectral method' for solving numerically such an eigenvalue problem. This
is an extension of the methods presented earlier in \cite{ach2008},
\cite{ach2009}.

\end{abstract}

\section{INTRODUCTION\label{sec1}}

We consider the numerical solution of \ the eigenvalue problem%
\begin{equation}
Lu(s)\equiv-\sum_{k,\ell=1}^{d}\frac{\partial}{\partial s_{k}}\left(
a_{k,\ell}(s)\frac{\partial u(s)}{\partial s_{\ell}}\right)  +\gamma
(s)u(s)=\lambda u(s),\quad\quad s\in\Omega\subseteq\mathbb{R}^{d}\label{e1}%
\end{equation}
with the Dirichlet boundary condition
\begin{equation}
u(s)\equiv0,\quad\quad s\in\partial\Omega.\label{e2}%
\end{equation}
Assume $d\geq2$. Let $\Omega$ be an open, simply--connected, and bounded
region in $\mathbb{R}^{d}$, and assume that its boundary $\partial\Omega$ is
smooth and sufficiently differentiable. Similarly, assume the functions
$\gamma(s)$ and $a_{i,j}(s)$, $1\leq i,j\leq d$, are several times
continuously differentiable over $\overline{\Omega}$. As usual, assume the
matrix $A(s)=\left[  a_{i,j}(s)\right]  $ is symmetric and satisfies the
strong ellipticity condition,%
\begin{equation}
\xi^{\text{T}}A(s)\xi\geq c_{0}\xi^{\text{T}}\xi,\quad\quad s\in
\overline{\Omega},\quad\xi\in\mathbb{R}^{d}\label{e2a}%
\end{equation}
with $c_{0}>0$. For convenience and without loss of generality, we assume
$\gamma(s)\geq0$, $s\in\Omega$; for otherwise, we can add a multiple of $u(s)$
to both sides of (\ref{e1}), shifting the eigenvalues by a known constant.

In the earlier papers \cite{ach2008} and \cite{ach2009} we introduced a
spectral method for the numerical solution of elliptic problems over $\Omega$
with Dirichlet and Neumann boundary conditions, respectively. In the present
work, this spectral method is extended to the numerical solution of the
eigenvalue problem for (\ref{e1})--(\ref{e2}), and in a later section it is
also extended to the Neumann problem%
\begin{align*}
-\Delta u(s)  &  =\lambda u(s),\quad\quad s\in\Omega\smallskip\\
\frac{\partial u}{\partial n}  &  =0,\quad\quad s\in\partial\Omega.
\end{align*}

\section{The Dirichlet problem\label{sec2}}

Our spectral method is based on polynomial approximation on the unit ball
$B_{d}$ in $\mathbb{R}^{d}$. To transform a problem defined on $\Omega$ to an
equivalent problem defined on $B_{d}$, we review some \ ideas from
\cite{ach2008} and \cite{ach2009}, modifying them as appropriate for this paper.

Assume the existence of a function%
\begin{equation}
\Phi:\overline{B}_{d}\underset{onto}{\overset{1-1}{\longrightarrow}}%
\overline{\Omega}\label{e3}%
\end{equation}
with $\Phi$ a twice--differentiable mapping, and let $\Psi=\Phi^{-1}%
:\overline{\Omega}\underset{onto}{\overset{1-1}{\longrightarrow}}\overline
{B}_{d}$. \ For $v\in L^{2}\left(  \Omega\right)  $, let%
\begin{equation}
\widetilde{v}(x)=v\left(  \Phi\left(  x\right)  \right)  ,\quad\quad
x\in\overline{B}_{d}\subseteq\mathbb{R}^{d}\label{e3a}%
\end{equation}
and conversely,%
\begin{equation}
v(s)=\widetilde{v}\left(  \Psi\left(  s\right)  \right)  ,\quad\quad
s\in\overline{\Omega}\subseteq\mathbb{R}^{d}.\label{e3b}%
\end{equation}
Assuming $v\in H^{1}\left(  \Omega\right)  $, we can show%
\[
\nabla_{x}\widetilde{v}\left(  x\right)  =J\left(  x\right)  ^{\text{T}}%
\nabla_{s}v\left(  s\right)  ,\quad\quad s=\Phi\left(  x\right)
\]
with $J\left(  x\right)  $ the Jacobian matrix for $\Phi$ over the unit ball
$B_{d}$,%
\begin{equation}
J(x)\equiv\left(  D\Phi\right)  (x)=\left[  \frac{\partial\varphi_{i}%
(x)}{\partial x_{j}}\right]  _{i,j=1}^{d},\quad\quad x\in\overline{B}%
_{d}.\label{e3c}%
\end{equation}
To use our method for problems over a region $\Omega$, it is necessary to know
explicitly the functions $\Phi$ and $J$. \ We assume%
\begin{equation}
\det J(x)\neq0,\quad\quad x\in\overline{B}_{d}.\label{e3d}%
\end{equation}
Similarly,%
\[
\nabla_{s}v(s)=K(s)^{\text{T}}\nabla_{x}\widetilde{v}(x),\quad\quad x=\Psi(s)
\]
with $K(s)$ the Jacobian matrix for $\Psi$ over $\Omega$. By differentiating
the identity
\[
\Psi\left(  \Phi\left(  x\right)  \right)  =x,\quad\quad x\in\overline{B}_{d}%
\]
we obtain%
\begin{equation}
K\left(  \Phi\left(  x\right)  \right)  =J\left(  x\right)  ^{-1}.\label{en16}%
\end{equation}
Assumptions about the differentiability of $\widetilde{v}\left(  x\right)  $
can be related back to assumptions on the differentiability of $v(s)$ and
$\Phi(x)$.

\begin{lemma}
If $\Phi\in C^{k}\left(  \overline{B}_{d}\right)  $ and $v\in C^{m}\left(
\overline{\Omega}\right)  $, then $\widetilde{v}\in C^{q}\left(  \overline
{B}_{d}\right)  $ with $q=\min\left\{  k,m\right\}  $.
\end{lemma}

\noindent\textbf{Proof}. A proof is straightforward using (\ref{e3a}).
$\ \left.  {}\right.  \hfill$\U{25a0}\medskip

\noindent A converse statement can be made as regards $\widetilde{v}$, $v$,
and $\Psi$ in (\ref{e3b}).

Consider now the nonhomogeneous problem $Lu=f$,%
\begin{equation}
Lu(s)\equiv-\sum_{k,\ell=1}^{d}\frac{\partial}{\partial s_{k}}\left(
a_{k,\ell}(s)\frac{\partial u(s)}{\partial s_{\ell}}\right)  +\gamma
(s)u(s)=f(s),\quad\quad s\in\Omega\subseteq\mathbb{R}^{d}.\label{e5}%
\end{equation}
Using the transformation (\ref{e3}), it is shown in \cite[Thm 2]{ach2008} that
(\ref{e5}) is equivalent to%
\begin{equation}%
\begin{array}
[c]{r}%
-%
{\displaystyle\sum\limits_{k,\ell=1}^{d}}
\dfrac{\partial}{\partial x_{k}}\left(  \widetilde{a}_{k,\ell}(x)\det\left(
J(x)\right)  \dfrac{\partial\widetilde{v}(x)}{\partial x_{\ell}}\right)
+\left[  \widetilde{\gamma}(x)\det J(x)\right]  \,\widetilde{u}(x)\smallskip\\
=\widetilde{f}\left(  x\right)  \det J(x),\quad\quad x\in B_{d}%
\end{array}
\label{e6}%
\end{equation}
with the matrix $\widetilde{A}\left(  x\right)  \equiv\left[  \widetilde{a}%
_{i,j}(x)\right]  $ given by%
\begin{equation}
\widetilde{A}\left(  x\right)  =J\left(  x\right)  ^{-1}A\left(  \Phi\left(
x\right)  \right)  J\left(  x\right)  ^{-\text{T}}.\label{e7}%
\end{equation}
The matrix $\widetilde{A}$ satisfies the analogue of (\ref{e2a}), but over
$B_{d}$. Thus the original eigenvalue problem (\ref{e1})--(\ref{e2}) can be
replaced by%
\begin{equation}%
\begin{array}
[c]{r}%
-%
{\displaystyle\sum\limits_{k,\ell=1}^{d}}
\dfrac{\partial}{\partial x_{k}}\left(  \widetilde{a}_{k,\ell}(x)\det\left(
J(x)\right)  \dfrac{\partial\widetilde{u}(x)}{\partial x_{\ell}}\right)
+\left[  \widetilde{\gamma}(x)\det J(x)\right]  \,\widetilde{u}(x)\smallskip\\
=\lambda\widetilde{u}(x)\det J(x),\quad\quad x\in B_{d}%
\end{array}
\label{e8}%
\end{equation}
As a consequence of this transformation, we can work with an elliptic problem
defined over $B_{d}$ rather than over the original region $\Omega$.

\subsection{The variational framework\label{sec2.1}}

To develop our numerical method, we need a variational framework for
(\ref{e5}) with the Dirichlet condition $u=0$ on $\partial\Omega$. \ As usual,
multiply both sides of (\ref{e5}) by an arbitary $v\in H_{0}^{1}\left(
\Omega\right)  $, integrate over $\Omega$, and apply integration by parts.
This yields the problem of finding $u\in H_{0}^{1}\left(  \Omega\right)  $
such that%
\begin{equation}
\mathcal{A}\left(  u,v\right)  =\left(  f,v\right)  \equiv\ell\left(
v\right)  ,\quad\quad\text{for all }v\in H_{0}^{1}\left(  \Omega\right)
\label{e8a}%
\end{equation}
with%
\begin{equation}
\mathcal{A}\left(  v,w\right)  =\int_{\Omega}\left[  \sum_{k,\ell=1}%
^{d}a_{k,\ell}(s)\frac{\partial v(s)}{\partial s_{\ell}}\frac{\partial
w(s)}{\partial s_{k}}+\gamma(s)v(s)w(s)\right]  ds,\quad v,w\in H_{0}%
^{1}\left(  \Omega\right)  .\label{en9}%
\end{equation}
The right side of (\ref{e8a}) uses the inner product $\left(  \cdot
,\cdot\right)  $ of $L^{2}\left(  \Omega\right)  $. The operators $L$ and
$\mathcal{A}$ are related by%
\begin{equation}
\left(  Lu,v\right)  =\mathcal{A}\left(  u,v\right)  ,\quad\quad u\in
H^{2}\left(  \Omega\right)  ,\quad v\in H_{0}^{1}\left(  \Omega\right)
,\label{e4}%
\end{equation}
an identity we use later. The function $\mathcal{A}$ is an inner product and
it satisfies%
\begin{equation}
\left\vert \mathcal{A}\left(  v,w\right)  \right\vert \leq c_{\mathcal{A}%
}\left\Vert v\right\Vert _{1}\left\Vert w\right\Vert _{1},\quad\quad v,w\in
H_{0}^{1}\left(  \Omega\right) \label{e4a}%
\end{equation}%
\begin{equation}
\mathcal{A}\left(  v,v\right)  \geq c_{e}\Vert v\Vert_{1}^{2},\quad\quad v\in
H_{0}^{1}\left(  \Omega\right) \label{e4b}%
\end{equation}
for some positive constants $c_{\mathcal{A}}$ and $c_{e}$.

Associated with the Dirichlet problem%
\begin{align}
Lu(x)  &  =f(x),\quad\quad x\in\Omega,\quad f\in L^{2}\left(  \Omega\right)
\label{e4c}\\
u(x)  &  =0,\quad\quad x\in\partial\Omega\label{e4d}%
\end{align}
is the Green's function integral operator%
\begin{equation}
u(x)=\mathcal{G}f(x).\label{e9}%
\end{equation}

\begin{lemma}
The operator $\mathcal{G}$ is a bounded and self--adjoint operator from
$L^{2}\left(  \Omega\right)  $ into $H_{0}^{2}\left(  \Omega\right)  $.
Moreover, it is a compact operator from $L^{2}\left(  \Omega\right)  $ into
$H_{0}^{1}\left(  \Omega\right)  $, and more particularly, it is a compact
operator from $H_{0}^{1}\left(  \Omega\right)  $ into $H_{0}^{1}\left(
\Omega\right)  $.
\end{lemma}

\noindent\textbf{Proof}. A proof can be based on \cite[\S 6.3, Thm. 5]{evans}
together with the fact that the embedding of $H_{0}^{2}\left(  \Omega\right)
$ into $H_{0}^{1}\left(  \Omega\right)  $ is compact. The symmetry follows
from the self--adjointness of the original problem (\ref{e4c})--(\ref{e4d}%
).$\ \left.  {}\right.  \hfill$\U{25a0}\medskip

We convert (\ref{e4}) to%
\begin{equation}
\left(  f,v\right)  =\mathcal{A}\left(  \mathcal{G}f,v\right)  ,\quad\quad
v\in H_{0}^{1}\left(  \Omega\right)  ,\quad f\in L^{2}\left(  \Omega\right)
.\label{e10}%
\end{equation}

The problem (\ref{e4c})--(\ref{e4d}) has the following variational
reformulation: find $u\in H_{0}^{1}\left(  \Omega\right)  $ such that
\begin{equation}
\mathcal{A}\left(  u,v\right)  =\ell(v),\quad\quad\forall v\in H_{0}%
^{1}\left(  \Omega\right)  .\label{e10a}%
\end{equation}
This problem can be shown to have a unique solution\ $u$ by using the
Lax--Milgram Theorem to imply its existence; see \cite[Thm. 8.3.4]%
{atkinson-han}. In addition,
\[
\Vert u\Vert_{1}\leq\frac{1}{c_{e}}\Vert\ell\Vert
\]
with $\Vert\ell\Vert$ denoting the operator norm for $\ell$ regarded as a
linear functional on $H_{0}^{1}\left(  \Omega\right)  $.

\subsection{The approximation scheme\label{sec2.2}}

Denote by $\Pi_{n}$ the space of polynomials in $d$ variables that are of
degree $\leq n$: $p\in\Pi_{n}$ if it has the form%
\[
p(x)=\sum_{\left\vert i\right\vert \leq n}a_{i}x_{1}^{i_{1}}x_{2}^{i_{2}}\dots
x_{d}^{i_{d}}%
\]
with $i$ a multi--integer, $i=\left(  i_{1},\dots,i_{d}\right)  $, and
$\left\vert i\right\vert =i_{1}+\cdots+i_{d}$. Over $B_{d}$, our approximation
subspace is%
\begin{equation}
\widetilde{\mathcal{X}}_{n}=\left\{  \left(  1-\left\Vert x\right\Vert
_{2}^{2}\right)  p(x)\mid p\in\Pi_{n}\right\} \label{e11}%
\end{equation}
with $\left\Vert x\right\Vert _{2}^{2}=x_{1}^{2}+\cdots+x_{d}^{2}$. The
subspaces $\Pi_{n}$ and $\widetilde{\mathcal{X}}_{n}$ have dimension%
\[
N\equiv N_{n}=\binom{n+d}{d}%
\]
However our problem (\ref{e8a}) is defined over $\Omega$, and thus we use a
modification of $\widetilde{\mathcal{X}}_{n}$:%
\begin{equation}
\mathcal{X}_{n}=\left\{  \psi\left(  s\right)  =\widetilde{\psi}\left(
\Psi\left(  s\right)  \right)  :\widetilde{\psi}\in\widetilde{\mathcal{X}}%
_{n}\right\} \label{e11a}%
\end{equation}
The finite dimensional set $\mathcal{X}_{n}\subseteq H_{0}^{1}\left(
\Omega\right)  $. This set of functions is used in the initial definition of
our numerical scheme and for its convergence analysis; but the simpler space
$\widetilde{\mathcal{X}}_{n}$ is used in the actual implementation of the method.

To solve (\ref{e10a}) (and thus (\ref{e4c})--(\ref{e4d})) approximately, we
use the Galerkin method with trial space $\mathcal{X}_{n}$ to find $u_{n}%
\in\mathcal{X}_{n}$ for which%
\begin{equation}
\mathcal{A}\left(  u_{n},v\right)  =\ell(v),\quad\quad\forall v\in
\mathcal{X}_{n}.\label{e23}%
\end{equation}
For the eigenvalue problem (\ref{e1}), find $u_{n}\in\mathcal{X}_{n}$ for
which
\begin{equation}
\mathcal{A}\left(  u_{n},v\right)  =\lambda\left(  u_{n},v\right)  ,\quad
\quad\forall v\in\mathcal{X}_{n}.\label{e24}%
\end{equation}
Write
\begin{equation}
u_{n}\left(  s\right)  =\sum_{j=1}^{N}\alpha_{j}\psi_{j}\left(  s\right)
\label{e25}%
\end{equation}
with $\left\{  \psi_{j}\right\}  _{j=1}^{N}$ a basis of $\mathcal{X}_{n}$.
Then (\ref{e24}) becomes%
\begin{equation}
\sum_{j=1}^{N}\alpha_{j}\mathcal{A}\left(  \psi_{j},\psi_{i}\right)
=\lambda\sum_{j=1}^{N}\alpha_{j}\left(  \psi_{j},\psi_{i}\right)  ,\quad
i=1,\dots,N\label{e26}%
\end{equation}

The coefficients can be related back to a polynomial basis for
$\widetilde{\mathcal{X}}_{n}$ and to integrals over $B_{d}$. Let $\left\{
\widetilde{\psi}_{j}\right\}  $ denote the basis of $\widetilde{\mathcal{X}%
}_{n}$ \ corresponding to the basis $\left\{  \psi_{j}\right\}  $\ for
$\mathcal{X}_{n}$. Using the transformation $s=\Phi(x)$,%
\begin{align}
\left(  \psi_{j},\psi_{i}\right)   &  =\int_{\Omega}\psi_{j}\left(  s\right)
\psi_{i}\left(  s\right)  \,ds\nonumber\\
&  =\int_{B_{d}}\widetilde{\psi}_{j}\left(  x\right)  \widetilde{\psi}%
_{i}\left(  x\right)  \left\vert \det J\left(  x\right)  \right\vert
\,dx\label{e26b}%
\end{align}%
\begin{align*}
\mathcal{A}\left(  \psi_{j},\psi_{i}\right)   &  =\int_{\Omega}\left[
\sum_{k,\ell=1}^{d}a_{k,\ell}\left(  s\right)  \frac{\partial\psi_{j}%
(s)}{\partial s_{k}}\frac{\partial\psi_{i}(s)}{\partial s_{\ell}}%
+\gamma(s)\psi_{j}(s)\psi_{i}(s)\right]  \,ds\\
&  =\int_{\Omega}\left[  \left\{  \nabla_{s}\psi_{i}\left(  s\right)
\right\}  ^{\text{T}}A(s)\left\{  \nabla_{s}\psi_{j}\left(  s\right)
\right\}  +\gamma(s)\psi_{j}(s)\psi_{i}(s)\right]  \,ds\\
&  =\int_{\Omega}\left[  \left\{  K(\Phi\left(  x\right)  )^{\text{T}}%
\nabla_{x}\widetilde{\psi}_{i}\left(  x\right)  \right\}  ^{\text{T}}A\left(
\Phi\left(  x\right)  \right)  \left\{  K(\Phi\left(  x\right)  )^{\text{T}%
}\nabla_{x}\widetilde{\psi}_{j}\left(  x\right)  \right\}  \right. \\
&  \left.  \quad\quad\quad+\widetilde{\gamma}(x)\widetilde{\psi}%
_{j}(x)\widetilde{\psi}_{i}(x)\right]  \left\vert \det J\left(  x\right)
\right\vert \,dx\\
&  =\int_{B_{d}}\left[  \nabla_{x}\widetilde{\psi}_{i}\left(  x\right)
^{\text{T}}\widetilde{A}(x)\nabla_{x}\widetilde{\psi}_{j}\left(  x\right)
+\widetilde{\gamma}(x)\widetilde{\psi}_{i}\left(  x\right)  \widetilde{\psi
}_{j}\left(  x\right)  \right]  \left\vert \det J\left(  x\right)  \right\vert
\,dx
\end{align*}
with the matrix $\widetilde{A}(x)$ given in (\ref{e7}). With these evaluations
of the coefficients, it is straightforward to show that (\ref{e26}) is
equivalent to a Galerkin method for (\ref{e7}) using the standard inner
product of $L^{2}\left(  B_{d}\right)  $ and the approximating subspace
$\widetilde{\mathcal{X}}_{n}$.

\subsection{Convergence analysis\label{sec2.3}}

The scheme (\ref{e26}) is implicitly a numerical approximation of the integral
equation eigenvalue problem
\begin{equation}
\lambda\mathcal{G}u=u.\label{e27}%
\end{equation}

\begin{lemma}
The numerical method (\ref{e24}) is equivalent to the Galerkin method
approximation of the integral equation (\ref{e27}), with the Galerkin method
based on the inner product $\mathcal{A}\left(  \cdot,\cdot\right)  $ for
$H_{0}^{1}\left(  \Omega\right)  $.
\end{lemma}

\noindent\textbf{Proof}. \ For the Galerkin solution of (\ref{e27}) we seek a
function $u_{n}$ in the form (\ref{e25}), and we force the residual to be
orthogonal to $\mathcal{X}_{n}$. This leads to
\begin{equation}
\lambda\sum_{j=1}^{N}\alpha_{j}\mathcal{A}\left(  \mathcal{G}\psi_{j},\psi
_{i}\right)  =\sum_{j=1}^{N}\alpha_{j}\mathcal{A}\left(  \psi_{j},\psi
_{i}\right) \label{e28}%
\end{equation}
for $i=1,\dots,N$. From (\ref{e10}), we have $\mathcal{A}\left(
\mathcal{G}\psi_{j},\psi_{i}\right)  =\left(  \psi_{j},\psi_{i}\right)  $, and
thus%
\[
\lambda\sum_{j=1}^{N}\alpha_{j}\left(  \psi_{j},\psi_{i}\right)  =\sum
_{j=1}^{N}\alpha_{j}\mathcal{A}\left(  \psi_{j},\psi_{i}\right)
\]
This is exactly the same as (\ref{e26}).$\ \left.  {}\right.  \hfill
$\U{25a0}\medskip

Let $\mathcal{P}_{n}$ be the orthogonal projection of $H_{0}^{1}\left(
B\right)  $ onto $\mathcal{X}_{n}$, based on the inner product $\mathcal{A}%
\left(  \cdot,\cdot\right)  $. Then (\ref{e28}) is the Galerkin approximation,%
\begin{equation}
\mathcal{P}_{n}\mathcal{G}u_{n}=\frac{1}{\lambda}u_{n},\quad u_{n}%
\in\mathcal{X}_{n}\label{e30}%
\end{equation}
for the integral equation eigenvalue problem (\ref{e27}). Much is known about
such schemes, as we discuss below. The conversion of the eigenvalue problem
(\ref{e24}) into the equivalent eigenvalue problem (\ref{e30}) is motivated by
a similar idea used in Osborn \cite{osborn}.

The numerical solution of eigenvalue problems for compact integral operators
has been studied by many people for over a century. With Galerkin methods, we
note particularly the early work of\ Krasnoselskii \cite[p. 178]%
{krasnoselskii}. The book of Chatelin \cite{chatelin} presents and summarizes
much of the literature on the numerical solution of such eigenvalue problems
for compact operators. For our work we use the results given in
\cite{atkinson67b}, \cite{atkinson75} for pointwise convergent\ operator
approximations that are collectively compact.

We begin with some preliminary lemmas.

\begin{lemma}
\label{Equiv1}For suitable positive constants $c_{1}$ and $c_{2}$,
\[
c_{1}\Vert\widetilde{v}\Vert_{H_{0}^{1}\left(  B_{d}\right)  }\leq\Vert
v\Vert_{H_{0}^{1}\left(  \Omega\right)  }\leq c_{2}\Vert\widetilde{v}%
\Vert_{H_{0}^{1}\left(  B_{d}\right)  }%
\]
for all functions $v\in H_{0}^{1}\left(  \Omega\right)  $, with $\widetilde{v}
$ the corresponding function of (\ref{e3a}). Thus, for a sequence $\left\{
v_{n}\right\}  $ in $H_{0}^{1}\left(  \Omega\right)  $,%
\begin{equation}
v_{n}\rightarrow v\text{\quad in }H_{0}^{1}\left(  \Omega\right)
\iff\widetilde{v}_{n}\rightarrow\widetilde{v}\text{\quad in }H_{0}^{1}\left(
B_{d}\right) \label{e17}%
\end{equation}
with $\left\{  \widetilde{v}_{n}\right\}  $ the corresponding sequence in
$H_{0}^{1}\left(  B_{d}\right)  $.
\end{lemma}

\noindent\textbf{Proof}. Begin by noting that there is a 1-1 correspondence
between $H_{0}^{1}\left(  \Omega\right)  $ and $H_{0}^{1}\left(  B_{d}\right)
$ based on using (\ref{e3})--(\ref{e3b}). Next,
\begin{align*}
\Vert v\Vert_{H_{0}^{1}\left(  \Omega\right)  }^{2}  &  =\int_{\Omega}\left[
\left\vert \nabla v\left(  s\right)  \right\vert ^{2}+\left\vert
v(s)\right\vert ^{2}\right]  ds\\
&  =\int_{B_{d}}\left[  \left\vert \nabla\widetilde{v}\left(  x\right)
^{\text{T}}J\left(  x\right)  ^{-1}J\left(  x\right)  ^{-\text{T}}%
\nabla\widetilde{v}\left(  x\right)  \right\vert +\left\vert \widetilde{v}%
(x)\right\vert ^{2}\right]  \left\vert \det J(x)\right\vert \,dx\\
&  \leq\left[  \max_{x\in B}\left\vert \det J(x)\right\vert \right]
\max\left\{  \max_{x\in B}\Vert J\left(  x\right)  ^{-1}\Vert^{2},1\right\}
\int_{B_{d}}\left[  \left\vert \nabla\widetilde{v}\left(  x\right)
\right\vert ^{2}+\left\vert \widetilde{v}(x)\right\vert ^{2}\right]  \,dx
\end{align*}%
\begin{equation}
\Vert v\Vert_{H_{0}^{1}\left(  \Omega\right)  }\leq c_{2}\Vert\widetilde{v}%
\Vert_{H_{0}^{1}\left(  B_{d}\right)  }\label{en16a}%
\end{equation}
for a suitable constant $c_{2}\left(  \Omega\right)  $. The reverse
inequality, with the roles of $\Vert\widetilde{v}\Vert_{H_{0}^{1}\left(
B_{d}\right)  }$ and $\Vert v\Vert_{H_{0}^{1}\left(  \Omega\right)  }$
reversed, follows by an analogous argument. $\ \left.  {}\right.  \hfill
$\U{25a0}\medskip

\begin{lemma}
The set $\cup_{n\geq1}\mathcal{X}_{n}$ is dense in $H_{0}^{1}\left(
\Omega\right)  $.
\end{lemma}

\noindent\textbf{Proof}. The set $\cup_{n\geq1}\widetilde{\mathcal{X}}_{n}$ is
dense in $H_{0}^{1}\left(  B_{d}\right)  $, a result shown in \cite[see
(15)]{ach2008}. We can then use the correspondence between $H_{0}^{1}\left(
\Omega\right)  $ and $H_{0}^{1}\left(  B_{d}\right)  $, given in Lemma
\ref{Equiv1},\ to show that $\cup_{n\geq1}\mathcal{X}_{n}$ is dense in
$H_{0}^{1}\left(  \Omega\right)  .\ \left.  {}\right.  \hfill$\U{25a0}\medskip

\begin{lemma}
The standard norm $\Vert\cdot\Vert_{1}$ on $H_{0}^{1}\left(  \Omega\right)  $
and the norm $\Vert v\Vert_{\mathcal{A}}=\sqrt{\mathcal{A}\left(  v,v\right)
} $ are equivalent in the topology they generate. More precisely,%
\begin{equation}
\sqrt{c_{e}}\Vert v\Vert_{1}\leq\Vert v\Vert_{\mathcal{A}}\leq\sqrt
{c_{\mathcal{A}}}\Vert v\Vert_{1},\quad\quad v\in H_{0}^{1}\left(
\Omega\right)  .\label{e18}%
\end{equation}
with the constants $c_{\mathcal{A}}$, $c_{e}$ taken from (\ref{e4a}) and
(\ref{e4b}), respectively. Convergence of sequences $\left\{  v_{n}\right\}  $
is equivalent in the two norms.
\end{lemma}

\noindent\textbf{Proof}. It is immediate from (\ref{e4b}) and (\ref{e4a}%
).$\ \left.  {}\right.  \hfill$\U{25a0}\medskip

\begin{lemma}
For the orthogonal projection operator $\mathcal{P}_{n}$,
\begin{equation}
\mathcal{P}_{n}v\rightarrow v\quad\quad\text{as\quad}n\rightarrow\infty
,\quad\quad\text{for all }v\in H_{0}^{1}\left(  \Omega\right)  .\label{e32}%
\end{equation}

\end{lemma}

\noindent\textbf{Proof}. This follows from the definition of an orthogonal
projection operator and using the result that $\cup_{n\geq1}\mathcal{X}_{n}$
is dense in $H_{0}^{1}\left(  \Omega\right)  $. $\ \left.  {}\right.  \hfill$\U{25a0}

\begin{corollary}
For the integral operator $\mathcal{G}$,%
\begin{equation}
\left\Vert \left(  I-\mathcal{P}_{n}\right)  \mathcal{G}\right\Vert
\rightarrow0\quad\text{as\quad}n\rightarrow\infty\label{e33}%
\end{equation}
using the norm for operators from $H_{0}^{1}\left(  \Omega\right)  $ into
$H_{0}^{1}\left(  \Omega\right)  $.
\end{corollary}

\noindent\textbf{Proof}. Consider $\mathcal{G}$ and $\mathcal{P}_{n}$ as
operators on $H_{0}^{1}\left(  \Omega\right)  $ into $H_{0}^{1}\left(
\Omega\right)  $. The result follows from the compactness of $\mathcal{G}$ and
the pointwise convergence in (\ref{e32}); \ see \cite[Lemma 3.1.2]%
{atkinson97}.$\ \left.  {}\right.  \hfill$\U{25a0}\medskip

\begin{lemma}
$\left\{  \mathcal{P}_{n}\mathcal{G}\right\}  $ is collectively compact on
$H_{0}^{1}\left(  \Omega\right)  $ .
\end{lemma}

\noindent\textbf{Proof}. This follows for all such families $\left\{
\mathcal{P}_{n}\mathcal{G}\right\}  $ with $\mathcal{G}$ compact on a Banach
space $\mathcal{Y}$ and $\left\{  \mathcal{P}_{n}\right\}  $ pointwise
convergent on $\mathcal{Y}$. To prove this requires showing
\[
\left\{  \mathcal{P}_{n}\mathcal{G}v\mid\Vert v\Vert_{1}\leq1,\,n\geq
1\right\}
\]
has compact closure in $H_{0}^{1}\left(  \Omega\right)  $. This can be done by
showing that the set is totally bounded. We omit the details of the
proof.$\ \left.  {}\right.  \hfill$\U{25a0}\medskip

Summarizing, $\left\{  \mathcal{P}_{n}\mathcal{G}\right\}  $ is a collectively
compact family that is pointwise convergent on $H_{0}^{1}\left(
\Omega\right)  $. With this, the results in \cite{atkinson67b},
\cite{atkinson75} can be applied to (\ref{e30}) as a numerical approximation
to the eigenvalue problem (\ref{e27}). We summarize the application of those
results to (\ref{e30}).

\begin{theorem}
\label{ThmEig1}Let $\lambda$ be an eigenvalue for the problem (\ref{e1}%
)--(\ref{e2}), say of multiplicity $\nu$, and let $\chi^{(1)},\dots,\chi
^{(\nu)}$ be a basis for the associated eigenfunction subspace. Let
$\varepsilon>0$ be chosen such that there are no other eigenvalues of
(\ref{e1})--(\ref{e2}) within a distance $\varepsilon$ of $\lambda$. Let
$\sigma_{n}$ denote the eigenvalue solutions of (\ref{e24}) that are within
$\varepsilon$ of $\lambda$. Then for all sufficiently large $n$, say $n\geq
n_{0}$, the sum of the multiplicities of the approximating eigenvalues within
$\sigma_{n}$ equals $\nu$. Moreover,
\begin{equation}
\max_{\lambda_{n}\in\sigma_{n}}\left\vert \lambda-\lambda_{n}\right\vert \leq
c\,\max_{1\leq k\leq\nu}\Vert\left(  I-\mathcal{P}_{n}\right)  \chi^{(k)}%
\Vert_{1}\label{e35}%
\end{equation}
Let $u$ be an eigenfunction of (\ref{e1})--(\ref{e2}) associated with
$\lambda$. Let $\mathcal{W}_{n}$ be the direct sum of the eigenfunction
subspaces associated with the eigenvalues $\lambda_{n}\in\sigma_{n}$, and let
$\left\{  u_{n}^{(1)},\dots,u_{n}^{(\nu)}\right\}  $ be a basis for
$\mathcal{W}_{n}$. Then there is a sequence
\[
u_{n}=\sum_{k=1}^{\nu}\alpha_{n,k}u_{n}^{(k)}\in\mathcal{W}_{n}%
\]
for which%
\begin{equation}
\left\Vert u-u_{n}\right\Vert _{1}\leq c\,\max_{1\leq k\leq\nu}\Vert\left(
I-\mathcal{P}_{n}\right)  \chi^{(k)}\Vert_{1}\label{e36}%
\end{equation}
for some constant $c>0$ dependent on $\lambda$.
\end{theorem}

\noindent\textbf{Proof}. This is a direct consequence of results in
\cite{atkinson67b}, \cite{atkinson75}, together with the compactness of
$\mathcal{G}$ on $H_{0}^{1}\left(  \Omega\right)  $. It also uses the
equivalence of norms given in (\ref{e18}). $\ \left.  {}\right.  \hfill
$\U{25a0}\medskip

The norms $\Vert\left(  I-\mathcal{P}_{n}\right)  \chi^{(k)}\Vert_{1}$ can be
bounded using results from Ragozin \cite{ragozin}, just as was done in
\cite{ach2008}. We begin with the following result from \cite{ragozin}. The
corresponding result that is needed with the Neumann problem can be obtained
from \cite{bbl}.

\begin{lemma}
\label{ragozin}Assume $w\in C^{k+2}\left(  \overline{B}_{d}\right)  $ for some
$k>0$, and assume $\left.  w\right\vert _{\partial B}=0$. \ Then there is a
polynomial $q_{n}\in\widetilde{\mathcal{X}}_{n}$ for which%
\begin{equation}
\left\Vert w-q_{n}\right\Vert _{\infty}\leq D\left(  k,d\right)  n^{-k}\left(
n^{-1}\left\Vert w\right\Vert _{\infty,k+2}+\omega\left(  w^{(k+2)}%
,1/n\right)  \right) \label{e48}%
\end{equation}

In this,
\[
\left\Vert w\right\Vert _{\infty,k+2}=\sum_{\left\vert i\right\vert \leq
k+2}\left\Vert \partial^{i}w\right\Vert _{\infty}%
\]%
\[
\omega\left(  g,\delta\right)  =\sup_{\left\vert x-y\right\vert \leq\delta
}\left\vert g\left(  x\right)  -g\left(  y\right)  \right\vert
\]%
\[
\omega\left(  w^{(k+2)},\delta\right)  =\sum_{\left\vert i\right\vert
=k+2}\omega\left(  \partial^{i}w,\delta\right)
\]

\end{lemma}

\begin{theorem}
\label{ThmEig2}Recall the notation and assumptions of Theorem \ref{ThmEig1}.
Assume the eigenfunction basis functions $\chi^{(k)}\in C^{m+2}\left(
\Omega\right)  $ and assume $\Phi\in C^{m+2}\left(  B_{d}\right)  $, for some
$m\geq1$. Then%
\[
\max_{\lambda_{n}\in\sigma_{n}}\left\vert \lambda-\lambda_{n}\right\vert
=\mathcal{O}\left(  n^{-m}\right)
\]%
\[
\left\Vert u-u_{n}\right\Vert _{1}=\mathcal{O}\left(  n^{-m}\right)
\]

\end{theorem}

\noindent\textbf{Proof}. Begin with (\ref{e35})--(\ref{e36}). To obtain the
bounds for $\Vert\left(  I-\mathcal{P}_{n}\right)  u^{(k)}\Vert_{1}$ given
above using Lemma \ref{ragozin}, refer to the argument given in \cite{ach2008}%
.$\ \left.  {}\right.  \hfill$\U{25a0}\medskip

\section{Implementation\label{sec3}}

Consider the implementation of the Galerkin method of (\ref{e24}) for the
eigenvalue problem (\ref{e1}). We are to find the function $u_{n}%
\in\mathcal{X}_{n}$ satisfying (\ref{e26}). To do so, we begin by selecting a
basis for $\Pi_{n}$ that is orthonormal in $L^{2}\left(  B_{d}\right)  $,
denoting it by $\left\{  \widetilde{\varphi}_{1},\dots,\widetilde{\varphi}%
_{N}\right\}  $, with $N\equiv N_{n}=\dim\Pi_{n}$. Choosing such an
orthonormal basis is an attempt to have the matrix associated with the left
side of the linear system in (\ref{e26}) be better conditioned. Next, let
\begin{equation}
\widetilde{\psi}_{i}(x)=\left(  1-\left\Vert x\right\Vert _{2}^{2}\right)
\widetilde{\varphi}_{i}(x),\quad\quad i=1,\dots,N_{n}\label{e75}%
\end{equation}
to form a basis for $\widetilde{\mathcal{X}}_{n}$. As in (\ref{e11a}), let
$\left\{  \psi_{1},\dots,\psi_{N}\right\}  $ be the corresponding basis of
$\mathcal{X}_{n}$.

We seek%
\begin{equation}
u_{n}(s)=\sum_{j=1}^{N}\alpha_{j}\psi_{j}(s)\label{e76}%
\end{equation}
Then following the change of variable $s=\Phi\left(  x\right)  $,
\ (\ref{e26}) becomes%
\begin{equation}%
\begin{array}
[c]{r}%
{\displaystyle\sum\limits_{j=1}^{N}}
\alpha_{j}%
{\displaystyle\int_{B_{d}}}
\left[  \nabla\widetilde{\psi}_{j}\left(  x\right)  ^{\text{T}}\widetilde{A}%
(x)\nabla\widetilde{\psi}_{i}\left(  x\right)  +\widetilde{\gamma
}(x)\widetilde{\psi}_{j}\left(  x\right)  \widetilde{\psi}_{i}\left(
x\right)  \right]  \left\vert \det J\left(  x\right)  \right\vert
\,dx\medskip\\
=\lambda%
{\displaystyle\sum\limits_{j=1}^{N}}
\alpha_{j}%
{\displaystyle\int_{B_{d}}}
\widetilde{\psi}_{j}\left(  x\right)  \widetilde{\psi}_{i}\left(  x\right)
\left\vert \det J\left(  x\right)  \right\vert \,dx,\quad\quad i=1,\dots,N
\end{array}
\label{e78}%
\end{equation}
We need to calculate the orthonormal polynomials and their first partial
derivatives; and we also need to approximate the integrals in the linear
system. For an introduction to the topic of multivariate orthogonal
polynomials, see Dunkl and Xu \cite{DX} and Xu \cite{xu2004}. For multivariate
quadrature over the unit ball in $\mathbb{R}^{d}$, see Stroud \cite{stroud}.

\subsection{The planar case\label{sec3.1}}

The dimension of $\Pi_{n}$ is
\begin{equation}
N_{n}=\frac{1}{2}\left(  n+1\right)  \left(  n+2\right) \label{e79}%
\end{equation}
For notation, we replace $x$ with $\left(  x,y\right)  $. How do we choose the
orthonormal basis $\left\{  \widetilde{\varphi}_{\ell}(x,y)\right\}  _{\ell
=1}^{N}$ for $\Pi_{n}$? Unlike the situation for the single variable case,
there are many possible orthonormal bases over $B=D$, the unit disk in
$\mathbb{R}^{2}$. We have chosen one that is particularly convenient for our
computations. These are the "ridge polynomials" introduced by Logan and\ Shepp
\cite{Loga} for solving an image reconstruction problem. We summarize here the
results needed for our work.

Let
\[
\mathcal{V}_{n}=\left\{  P\in\Pi_{n}:\left(  P,Q\right)  =0\quad\forall
Q\in\Pi_{n-1}\right\}
\]
the polynomials of degree $n$ that are orthogonal to all elements of
$\Pi_{n-1}$. Then the dimension of $\mathcal{V}_{n}$ is $n+1$; moreover,%
\begin{equation}
\Pi_{n}=\mathcal{V}_{0}\oplus\mathcal{V}_{1}\oplus\cdots\oplus\mathcal{V}%
_{n}\label{e100}%
\end{equation}
It is standard to construct orthonormal bases of each $\mathcal{V}_{n}$ and to
then combine them to form an orthonormal basis of $\Pi_{n}$ using the latter
decomposition. \ As an orthonormal basis of $\mathcal{V}_{n}$ we use%
\begin{equation}
\widetilde{\varphi}_{n,k}(x,y)=\frac{1}{\sqrt{\pi}}U_{n}\left(  x\cos\left(
kh\right)  +y\sin\left(  kh\right)  \right)  ,\quad\left(  x,y\right)  \in
D,\quad h=\frac{\pi}{n+1}\label{e101}%
\end{equation}
for $k=0,1,\dots,n$. The function $U_{n}$ is the Chebyshev polynomial of the
second kind of degree $n$:%
\begin{equation}
U_{n}(t)=\frac{\sin\left(  n+1\right)  \theta}{\sin\theta},\quad\quad
t=\cos\theta,\quad-1\leq t\leq1,\quad n=0,1,\dots\label{e102}%
\end{equation}
The family $\left\{  \widetilde{\varphi}_{n,k}\right\}  _{k=0}^{n}$ is an
orthonormal basis of $\mathcal{V}_{n}$. As a basis of $\Pi_{n}$, we order
$\left\{  \widetilde{\varphi}_{n,k}\right\}  $ lexicographically based on the
ordering in (\ref{e101}) and (\ref{e100}):%
\[
\left\{  \widetilde{\varphi}_{\ell}\right\}  _{\ell=1}^{N}=\left\{
\widetilde{\varphi}_{0,0},\,\widetilde{\varphi}_{1,0},\,\widetilde{\varphi
}_{1,1},\,\widetilde{\varphi}_{2,0},\,\dots,\,\widetilde{\varphi}%
_{n,0},\,\dots,\widetilde{\varphi}_{n,n}\right\}
\]

Returning to (\ref{e75}), we define%
\begin{equation}
\widetilde{\psi}_{n,k}(x,y)=\left(  1-x^{2}-y^{2}\right)  \widetilde{\varphi
}_{n,k}(x,y)\label{e103}%
\end{equation}
To calculate the first order partial derivatives of $\widetilde{\psi}%
_{n,k}(x,y)$, we need $U_{n}^{^{\prime}}(t)$. The values of $U_{n}(t)$ and
$U_{n}^{^{\prime}}(t)$ are evaluated using the standard triple recursion
relations%
\begin{align*}
U_{n+1}(t)  &  =2tU_{n}(t)-U_{n-1}(t)\\
U_{n+1}^{^{\prime}}(t)  &  =2U_{n}(t)+2tU_{n}^{^{\prime}}(t)-U_{n-1}%
^{^{\prime}}(t)
\end{align*}

For the numerical approximation of the integrals in (\ref{e78}), which are
over $B$ being the unit disk, we use the formula%
\begin{equation}
\int_{B}g(x,y)\,dx\,dy\approx\sum_{l=0}^{q}\sum_{m=0}^{2q}g\left(  r_{l}%
,\frac{2\pi\,m}{2q+1}\right)  \omega_{l}\frac{2\pi}{2q+1}r_{l}\label{e106}%
\end{equation}
Here the numbers $\omega_{l}$ are the weights of the $\left(  q+1\right)
$-point Gauss-Legendre quadrature formula on $[0,1]$. Note that
\[
\int_{0}^{1}p(x)dx=\sum_{l=0}^{q}p(r_{l})\omega_{l},
\]
for all single-variable polynomials $p(x)$ with $\deg\left(  p\right)
\leq2q+1 $. The formula (\ref{e106}) uses the trapezoidal rule with $2q+1$
subdivisions for the integration over $\overline{B}$ in the azimuthal
variable. This quadrature (\ref{e106}) is exact for all polynomials $g\in
\Pi_{2q}$. This formula is also the basis of the hyperinterpolation formula
discussed in \cite{hac}.

\subsection{The three--dimensional case\label{sec3.2}}

In $\mathbb{R}^{3}$, the dimension of $\Pi_{n}$ is
\[
N_{n}=\binom{n+3}{3}=\frac{1}{6}\left(  n+1\right)  \left(  n+2\right)
\left(  n+3\right)
\]
Here we choose orthonormal polynomials on the unit ball as described in
\cite{DX},
\begin{align}
\widetilde{\varphi}_{m,j,\beta}(x)  &  =c_{m,j}p_{j}^{(0,m-2j+\frac{1}{2}%
)}(2\Vert x\Vert^{2}-1)S_{\beta,m-2j}\left(  \frac{x}{\Vert x\Vert}\right)
\medskip\nonumber\\
&  =c_{m,j}\Vert x\Vert^{m-2j}p_{j}^{(0,m-2j+\frac{1}{2})}(2\Vert x\Vert
^{2}-1)S_{\beta,m-2j}\left(  \frac{x}{\Vert x\Vert}\right)  ,\medskip
\label{e1001}\\
j  &  =0,\ldots,\lfloor m/2\rfloor,\quad\beta=0,1,\ldots,2(m-2j),\quad
m=0,1,\ldots,n\nonumber
\end{align}
Here $c_{m,j}=2^{\frac{5}{4}+\frac{m}{2}-j}$ is a constant, and $p_{j}%
^{(0,m-2j+\frac{1}{2})}$, $j\in\mathbb{N}_{0}$, are the normalized Jabobi
polynomials which are orthonormal on $[-1,1]$ with respect to the inner
product
\[
(v,w)=\int_{-1}^{1}(1+t)^{m-2j+\frac{1}{2}}v(t)w(t)\;dt,
\]
see for example \cite{abramowitz}, \cite{gautschi}. The functions
$S_{\beta,m-2j}$ are spherical harmonic functions, and they are given in
spherical coordinates by
\[
S_{\beta,k}(\phi,\theta)=\widetilde{c}_{\beta,k}\left\{
\begin{array}
[c]{ll}%
\cos(\frac{\beta}{2}\phi)T_{k}^{\frac{\beta}{2}}(\cos\theta),\medskip &
\beta\text{ even}\\
\sin(\frac{\beta+1}{2}\phi)T_{k}^{\frac{\beta+1}{2}}(\cos\theta), &
\beta\text{ odd}%
\end{array}
\right.
\]
The constant $\widetilde{c}_{\beta,k}$ is chosen in such a way that the
functions are orthonormal on the unit sphere $S^{2}$ in $\mathbb{R}^{3}$:
\[
\int_{S^{2}}S_{\beta,k}(x)\,S_{\widetilde{\beta},\widetilde{k}}(x)\,dS=\delta
_{\beta,\widetilde{\beta}}\,\delta_{k,\widetilde{k}}%
\]
The functions $T_{k}^{l}$ are the associated Legendre polynomials, see
\cite{hobson}, \cite{macroberts}. According to (\ref{e75}) we define the basis
for our space of trial functions by
\[
\widetilde{\psi}_{m,j,\beta}(x)=(1-\Vert x\Vert^{2})\widetilde{\varphi
}_{m,j,\beta}(x)
\]
and we can order the basis lexicographically. To calculate all of the above
functions we can use recursive algorithms similar to the one used for the
Chebyshev polynomials. These algorithms also allow the calculation of the
derivatives of each of these functions, see \cite{gautschi}, \cite{zhang-jin}

For the numerical approximation of the integrals in (\ref{e78}) we use a
quadrature formula for the unit ball $B$
\begin{align}
\int_{B}g(x)\,dx  &  =\int_{0}^{1}\int_{0}^{2\pi}\int_{0}^{\pi}\widetilde{g}%
(r,\theta,\phi)\,r^{2}\sin(\phi)\,d\phi\,d\theta\,dr\approx Q_{q}%
[g]\medskip\nonumber\\
Q_{q}[g]  &  :=\sum_{i=1}^{2q}\sum_{j=1}^{q}\sum_{k=1}^{q}\frac{\pi}%
{q}\,\omega_{j}\,\nu_{k}\widetilde{g}\left(  \frac{\zeta_{k}+1}{2},\frac
{\pi\;i}{2q},\arccos(\xi_{j})\right) \label{e1004}%
\end{align}
Here $\widetilde{g}(r,\theta,\phi)=g(x)$ is the representation of $g$ in
spherical coordinates. For the $\theta$ integration we use the trapezoidal
rule, because the function is $2\pi-$periodic in $\theta$. For the $r$
direction we use the transformation
\begin{align*}
\int_{0}^{1}r^{2}v(r)\;dr  &  =\int_{-1}^{1}\left(  \frac{t+1}{2}\right)
^{2}v\left(  \frac{t+1}{2}\right)  \frac{dt}{2}\medskip\\
&  =\frac{1}{8}\int_{-1}^{1}(t+1)^{2}v\left(  \frac{t+1}{2}\right)
\;dt\medskip\\
&  \approx\sum_{k=1}^{q}\underset{_{=:\nu_{k}}}{\underbrace{\frac{1}{8}\nu
_{k}^{\prime}}}v\left(  \frac{\zeta_{k}+1}{2}\right)
\end{align*}
where the $\nu_{k}^{\prime}$ and $\zeta_{k}$ are the weights and the nodes of
the Gauss quadrature with $q$ nodes on $[-1,1]$ with respect to the inner
product
\[
(v,w)=\int_{-1}^{1}(1+t)^{2}v(t)w(t)\,dt
\]
The weights and nodes also depend on $q$ but we omit this index. For the
$\phi$ direction we use the transformation
\begin{align*}
\int_{0}^{\pi}\sin(\phi)v(\phi)\,d\phi &  =\int_{-1}^{1}v(\arccos
(\phi))\,d\phi\medskip\\
&  \approx\sum_{j=1}^{q}\omega_{j}v(\arccos(\xi_{j}))
\end{align*}
where the $\omega_{j}$ and $\xi_{j}$ are the nodes and weights for the
Gauss--Legendre quadrature on $[-1,1]$. For more information on this
quadrature rule on the unit ball in $\mathbb{R}^{3}$, see \cite{stroud}.

Finally we need the gradient in Cartesian coordinates to approximate the
integral in (\ref{e78}), but the function $\widetilde{\varphi}_{m,j,\beta}(x)
$ in (\ref{e1001}) is given in spherical coordinates. Here we simply use the
chain rule, with $x=\left(  x,y,z\right)  $,
\begin{align*}
\frac{\partial}{\partial x}v(r,\theta,\phi)  &  =\frac{\partial}{\partial
r}v(r,\theta,\phi)\cos(\theta)\sin(\phi)-\frac{\partial}{\partial\theta
}v(r,\theta,\phi)\frac{\sin(\theta)}{r\sin(\phi)}\medskip\\
&  +\frac{\partial}{\partial\phi}v(r,\theta,\phi)\frac{\cos(\theta)\cos(\phi
)}{r}%
\end{align*}
and similarly for $\frac{\partial}{\partial y}$ and $\frac{\partial}{\partial
z}$.

\section{Numerical example\label{num_example}}%

\begin{figure}[tb]%
\centering
\includegraphics[
height=3in,
width=3.9998in
]%
{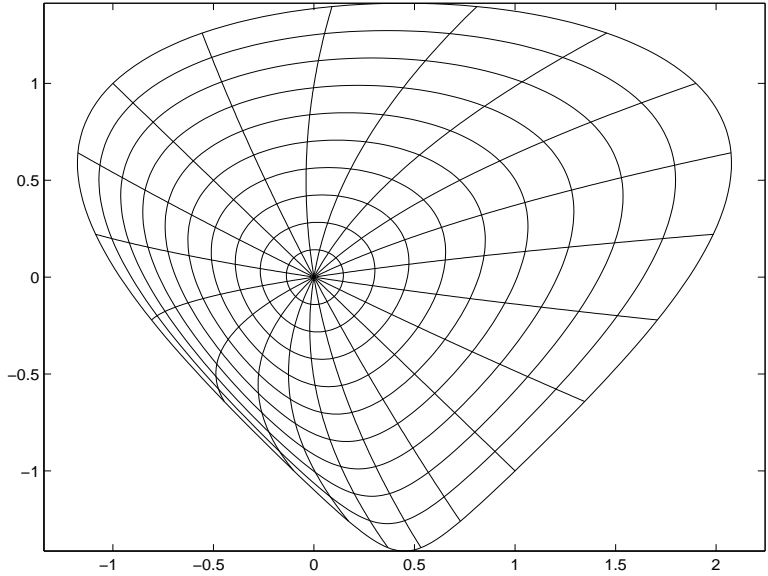}%
\caption{Images of (\ref{e132}), with $a=0.5$, for lines of constant radius
and constant azimuth on the unit disk.}%
\label{2D_Region}%
\end{figure}

Our programs are written in \textsc{Matlab}. The transformations have been so
chosen that we can invert explicitly the mapping $\Phi$, to be able to better
construct our test examples. This is not needed when applying the method; but
it simplified the construction of our test cases. The eigenvalue problem being
solved is
\begin{equation}
Lu(\mathbf{s})\equiv-\Delta u=\lambda u(\mathbf{s}),\quad\quad\mathbf{s}%
\in\text{$\Omega$}\subseteq\mathbb{R}^{d}\label{e130}%
\end{equation}
which corresponds to choosing $A=I$. Then we need to calculate%
\begin{equation}
\widetilde{A}\left(  \mathbf{x}\right)  =J\left(  \mathbf{x}\right)
^{-1}J\left(  \mathbf{x}\right)  ^{-\text{T}}\label{e131}%
\end{equation}

\subsection{The planar case\label{sec4.1}}%

\begin{figure}[tb]%
\centering
\includegraphics[
height=3in,
width=3.9998in
]%
{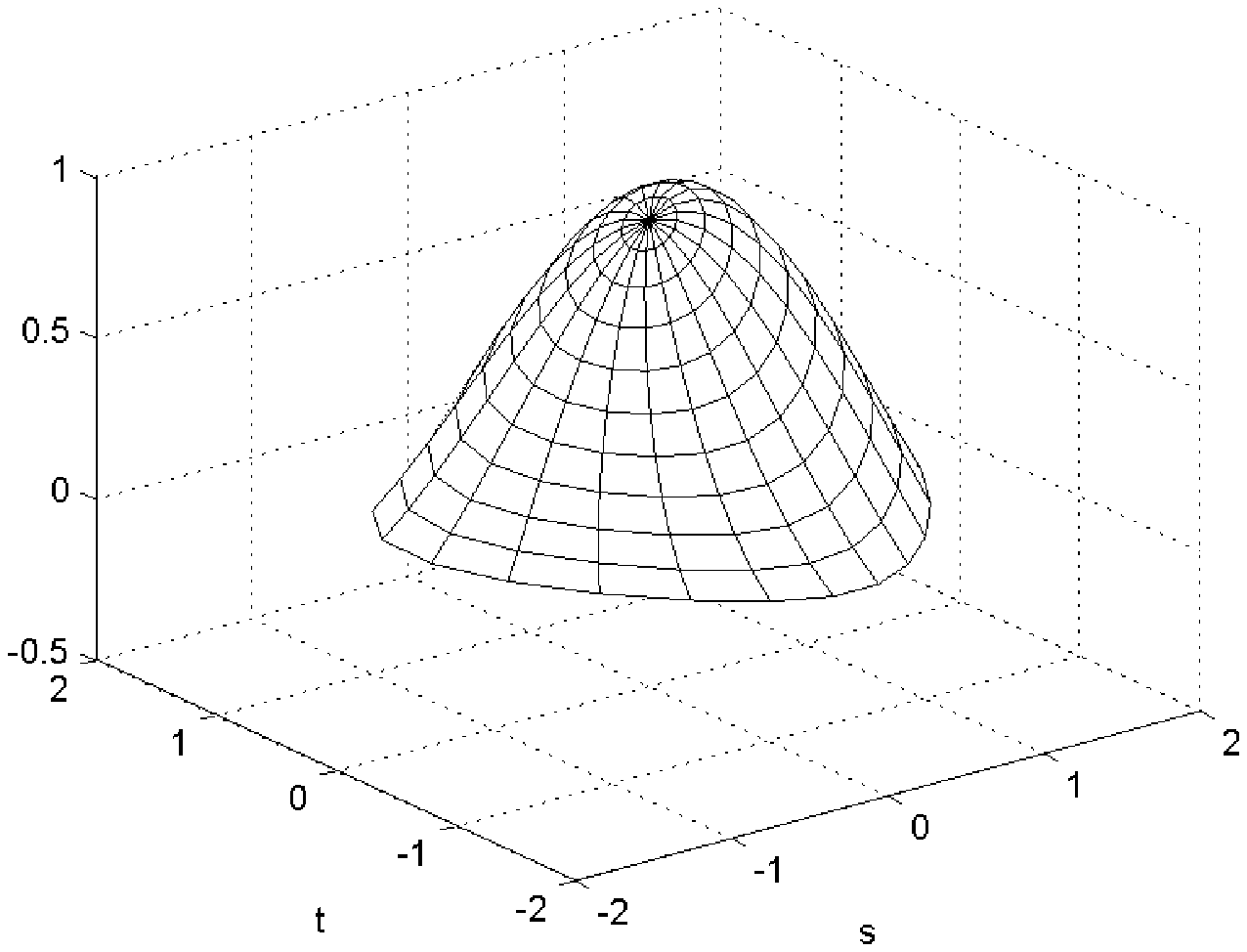}%
\caption{Eigenfunction corresponding to the approximate eigenvalue
\newline$\lambda_{1}\doteq2.96185.$}%
\label{EF1}%
\end{figure}

For our variables, we replace a point $x\in B_{2}$ with $\left(  x,y\right)
$, and we replace a point $s\in\Omega$ with $\left(  s,t\right)  $. \ Define
the mapping $\Phi:\overline{B}_{2}\rightarrow\overline{\Omega}$ by $\left(
s,t\right)  =\Phi\left(  x,y\right)  $,%
\begin{equation}%
\begin{array}
[c]{l}%
s=x-y+ax^{2}\\
t=x+y
\end{array}
\label{e132}%
\end{equation}
with $0<a<1$. It can be shown that $\Phi$ is a 1-1 mapping from the unit disk
$\overline{B}$. In particular, the inverse mapping $\Psi:\overline{\Omega
}\rightarrow\overline{B}$ is given by%
\begin{equation}%
\begin{array}
[c]{l}%
x=\dfrac{1}{a}\left[  -1+\sqrt{1+a\left(  s+t\right)  }\right]  \medskip\\
y=\dfrac{1}{a}\left[  at-\left(  -1+\sqrt{1+a\left(  s+t\right)  }\right)
\right]
\end{array}
\label{e133}%
\end{equation}
In Figure \ref{2D_Region}, we give the images in $\overline{\Omega}$ of the
circles $r=j/10$, $j=1,\dots,10$ and the azimuthal lines $\theta=j\pi/10$,
$j=1,\dots,20$.

The following information is needed when implementing the transformation from
$-\Delta u=\lambda u$ on $\Omega$ to a new equation on $B_{2}$:%
\[
D\Phi=J\left(  x,y\right)  =\left(
\begin{array}
[c]{cc}%
1+2ax & -1\\
1 & 1
\end{array}
\right)
\]%
\[
\det\left(  J\right)  =2\left(  1+ax\right)
\]%
\[
J\left(  x\right)  ^{-1}=\frac{1}{2\left(  1+ax\right)  }\left(
\begin{array}
[c]{cc}%
1 & 1\\
-1 & 1+2ax
\end{array}
\right)
\]%
\[
A=J\left(  x\right)  ^{-1}J\left(  x\right)  ^{-\text{T}}=\frac{1}{2\left(
1+ax\right)  ^{2}}\left(
\begin{array}
[c]{cc}%
1 & ax\\
ax & 2a^{2}x^{2}+2ax+1
\end{array}
\right)
\]

\begin{figure}[tb]%
\centering
\includegraphics[
height=3in,
width=3.9998in
]%
{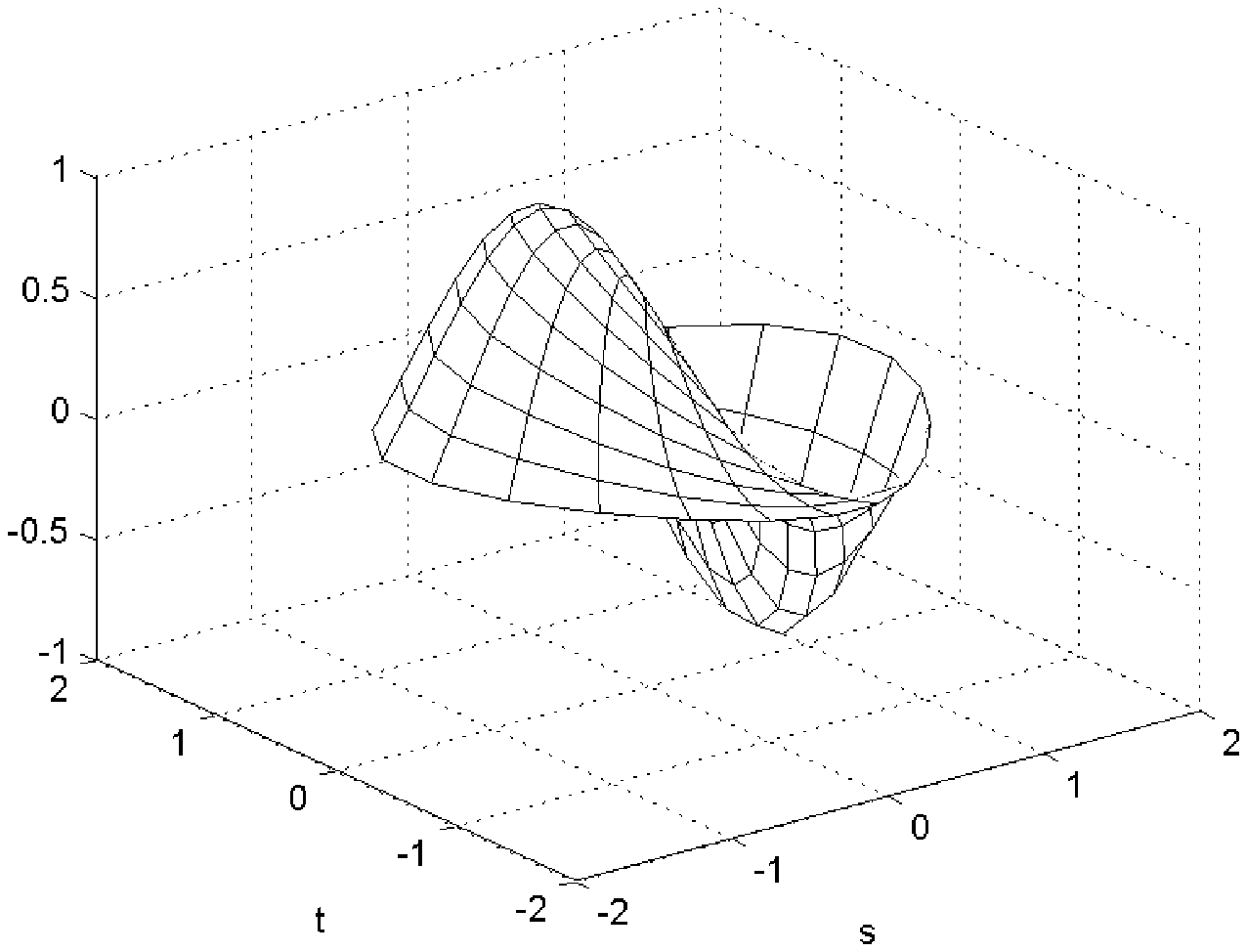}%
\caption{Eigenfunction corresponding to the approximate eigenvalue
\newline$\lambda_{2}\doteq7.24761.$}%
\label{EF2}%
\end{figure}

We give an example for this region $\Omega\,\ $with $a=0.5$. \ Figures
\ref{EF1} and \ref{EF2} contain the computed eigenfunctions for the two
smallest eigenvalues; these are based on the degree $n=8$ approximation.%

\begin{figure}[tb]%
\centering
\includegraphics[
height=3in,
width=3.9998in
]%
{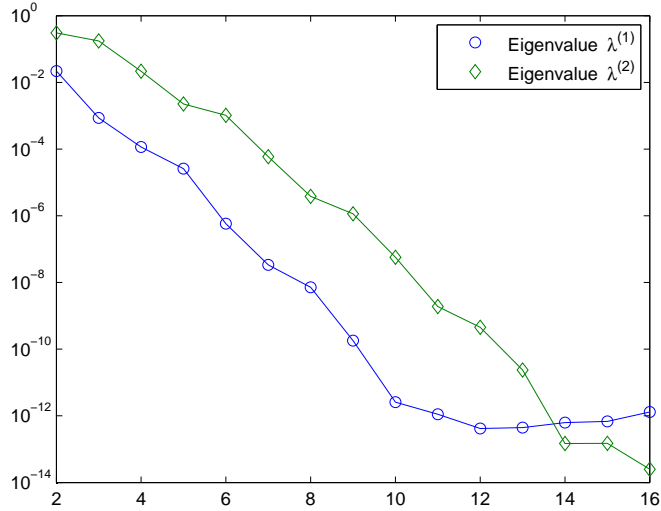}%
\caption{The values of $\left\vert \lambda_{n+1}^{(k)}-\lambda_{n}%
^{(k)}\right\vert $ for $k=1,2$ for increasing degree $n$. }%
\label{EV_Diff}%
\end{figure}

\begin{figure}[tb]%
\centering
\includegraphics[
height=3in,
width=3.9998in
]%
{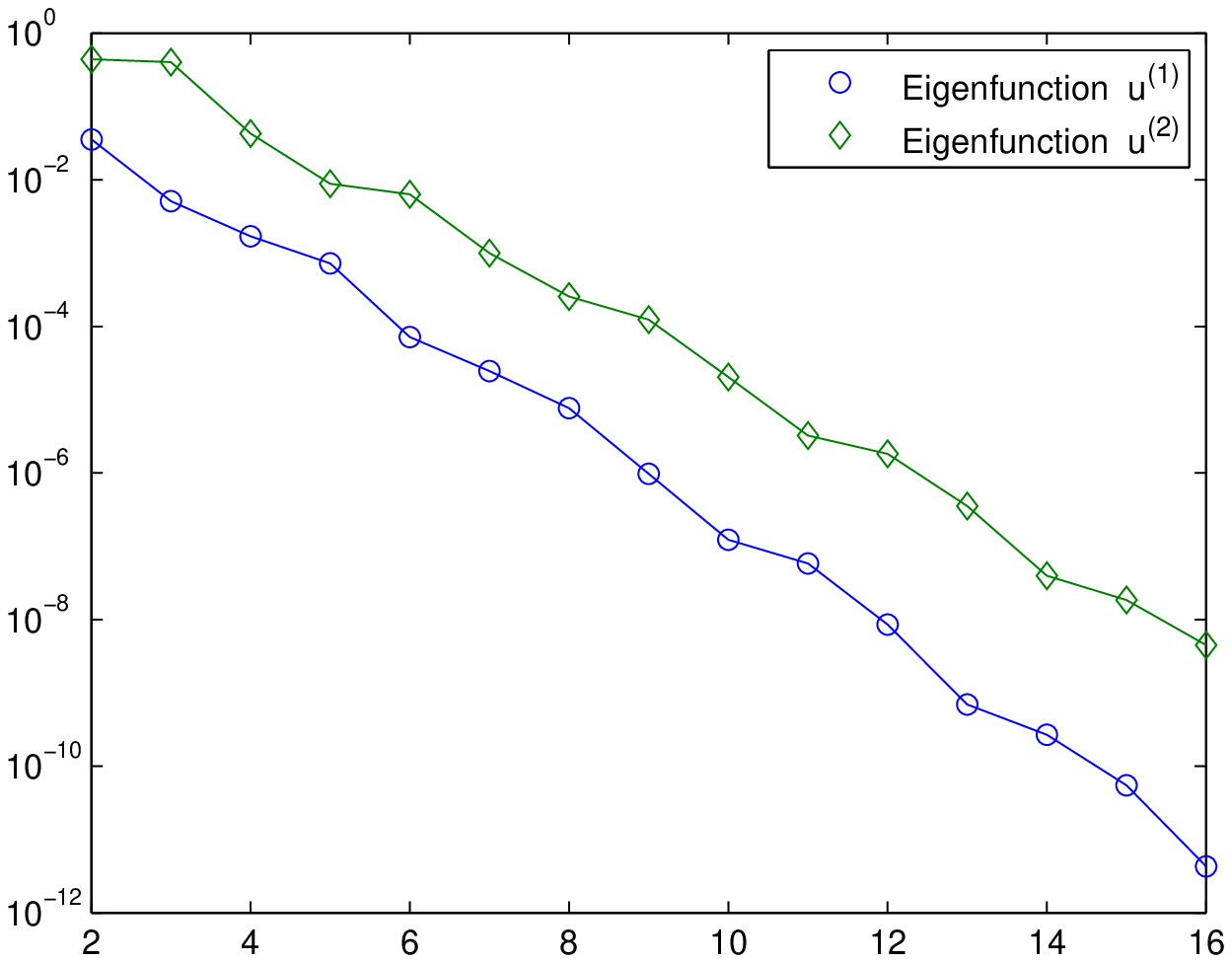}%
\caption{The values of $\Vert u_{n+1}^{(k)}-u_{n}^{(k)}\Vert_{\infty}$ for
$k=1,2$ for increasing degree $n$. }%
\label{EF_Diff}%
\end{figure}

\begin{figure}[tb]%
\centering
\includegraphics[
height=3in,
width=3.9998in
]%
{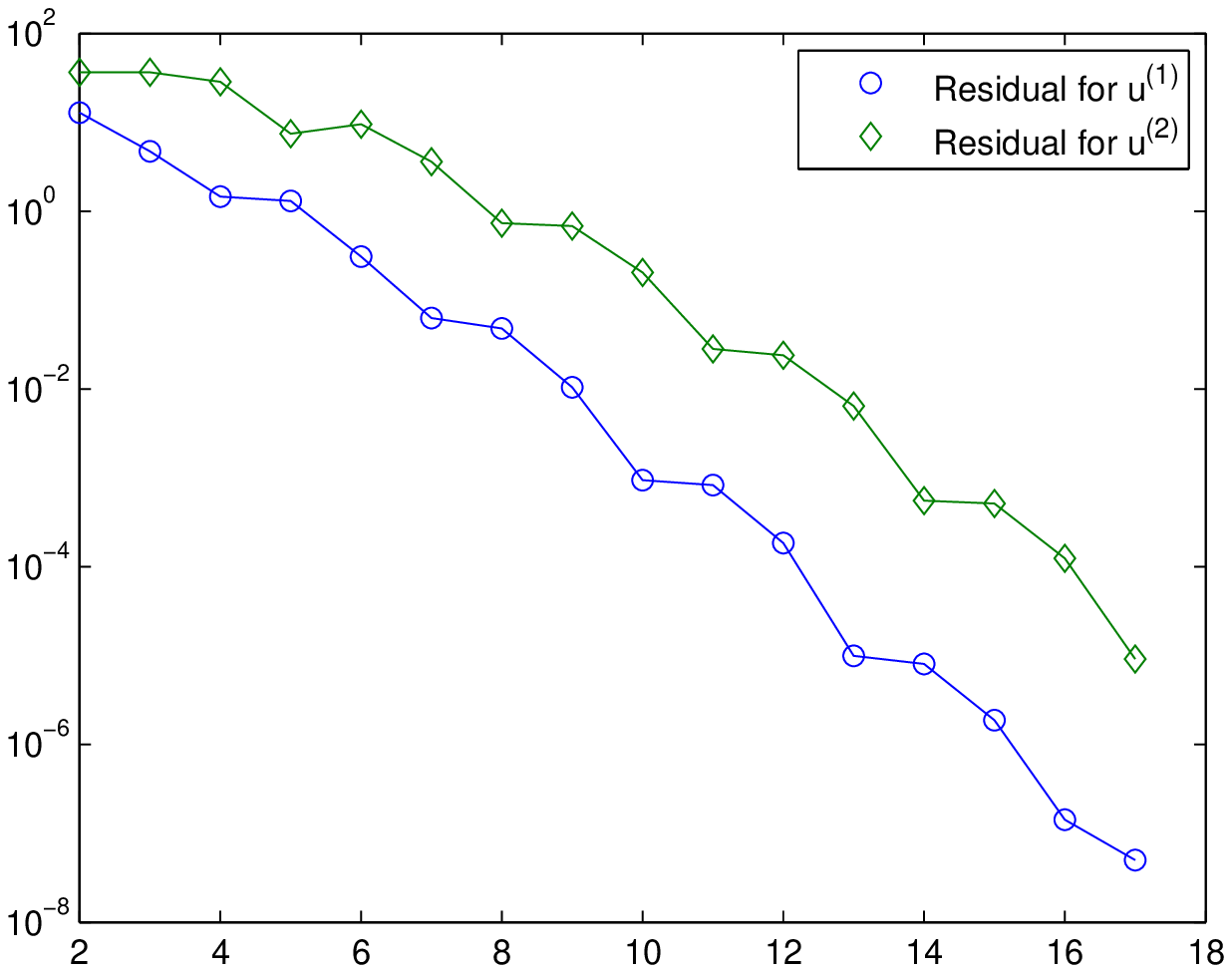}%
\caption{The values of $\Vert R_{n}^{(k)}\Vert_{\infty}$ for $k=1,2$ for
increasing degree $n$. }%
\label{Resid}%
\end{figure}

Because the true eigenfunctions and eigenvalues are unknown for almost all
cases (with the unit ball as an exception), we used other methods for studying
experimentally the rate of convergence. Let $\lambda_{n}^{(k)}$ denote the
value of the $k^{th}$ eigenvalue based on the degree $n$ polynomial
approximation, with the eigenvalues taken in increasing order. Let
$u_{n}^{(k)}$ denote a corresponding eigenfunction,
\[
\widetilde{u}_{n}^{(k)}\left(  x\right)  =\sum_{j=1}^{N_{n}}\alpha_{j}%
^{(n)}\widetilde{\psi}_{j}(x)
\]
with $\alpha^{(n)}\equiv\left[  \alpha_{1}^{(n)},\dots,\alpha_{N}%
^{(n)}\right]  $\ the eigenvector of (\ref{e78}) associated with the
eigenvalue \ $\lambda_{n}^{(k)}$. We normalize the eigenvectors by requiring
$\Vert\alpha^{(n)}\Vert_{\infty}=1$. Define%
\[
\Lambda_{n}=\left\vert \lambda_{n+1}^{(k)}-\lambda_{n}^{(k)}\right\vert
\]%
\[
D_{n}=\Vert u_{n+1}^{(k)}-u_{n}^{(k)}\Vert_{\infty}%
\]
Figures \ref{EV_Diff} and \ref{EF_Diff} \ show the decrease, respectively, of
$\Lambda_{n}$ and \ $D_{n}$ as $n$ increases. In both cases, we use a semi-log
scale. Also, consider the residual%
\[
R_{n}^{(k)}=-\Delta u_{n}^{(k)}-\lambda_{n}^{(k)}u_{n}^{(k)}%
\]
Figure \ \ref{Resid} shows the decrease of $\Vert R_{n}^{(k)}\Vert_{\infty} $,
again on a semi-log scale.%

%

%

These numerical results all indicate an exponential rate of convergence as a
function of the degree $n$ of the approximations $\left\{  \lambda_{n}%
^{(k)}:n\geq1\right\}  $ and $\left\{  u_{n}^{(k)}:n\geq1\right\}  $. In
Figure \ref{EV_Diff}, the maximum accuracy for $\lambda^{(1)}$ appears to have
been found with the degree $n=12$, approximately. For larger degrees, rounding
errors dominate. We also see that the accuracy for the first
eigenvalue-eigenfunction pair is better than that for the second such pair.

\subsection{The three--dimensional case\label{sec4.2}}

Here we consider the problem of finding eigenvalues and eigenfunctions for the
Neumann problem in $\Omega\subset{\mathbb{R}}^{3}$:
\begin{equation}%
\begin{array}
[c]{rcll}%
-\Delta u(s) & = & \lambda u(s), & s\in\Omega\smallskip\\
\dfrac{\partial u(s)}{\partial n} & = & 0, & s\in\partial\Omega
\end{array}
\label{eq1000}%
\end{equation}
Problem (\ref{eq1000}) is equivalent to
\begin{equation}%
\begin{array}
[c]{rcll}%
-\Delta u(s)+u(s) & = & (\lambda+1)u(s), & s\in\Omega\smallskip\\
\dfrac{\partial u(s)}{\partial n} & = & 0, & s\in\partial\Omega
\end{array}
\label{eq1001}%
\end{equation}
and $-\Delta+I:\,D_{N}\mapsto L^{2}(\Omega)$ is an invertible self--adjoint
operator with
\[
D_{L}=\left\{  u\in H^{2}(\Omega)\,\left\vert \,\frac{\partial u(s)}{\partial
n}=0,\,s\in\partial\Omega\right.  \right\}
\]
So there is a continuous solution operator $G:L^{2}(\Omega)\mapsto D_{L}$,
such that
\[
(-\Delta+I)\circ G|_{L^{2}(\Omega)}=I
\]
with $I$ the identity operator on $L^{2}\left(  \Omega\right)  $. If we
consider $G:H^{1}(\Omega)\mapsto H^{1}(\Omega)$, then $G$ is a compact
operator, because of the compact imbedding $H^{1}(\Omega)\hookrightarrow
L^{2}(\Omega)$ or $H^{2}(\Omega)\hookrightarrow H^{1}(\Omega)$; see
\cite{LadUra} or \cite{Triebel}.

We follow now Section \ref{sec2.1} to present the variational framework. A
solution of the inhomogeneous problem
\begin{equation}
Lu=f,\quad\quad f\in L^{2}(\Omega)\label{eq1002}%
\end{equation}
satisfies
\[
\int_{\Omega}\left(  -\sum_{j=1}^{3}\frac{\partial^{2}}{\partial s_{j}^{2}%
}u(s)+u(s)\right)  v(s)\,ds=\int_{\Omega}f(s)v(s)\,ds\quad\text{for all }v\in
H^{1}(\Omega)
\]
Applying integration by parts and using the fact that the normal derivative of
$u\in D_{L}$ is zero on $\partial\Omega$ we derive
\[
\int_{\Omega}\nabla_{s}u(s)\nabla_{s}v(s)+u(s)v(s)\,ds=\int_{\Omega
}f(s)v(s)\,ds\quad\text{for all }v\in H^{1}(\Omega)
\]
We denote the left hand side of this equation by $\mathcal{A}(u,v)$ and from
the Cauchy--Schwartz inequality we derive
\[
\mathcal{A}(u,v)\leq\Vert u\Vert_{H^{1}(\Omega)}\,\Vert v\Vert_{H^{1}(\Omega)}%
\]
and we have the equality
\[
\mathcal{A}(u,u)=\Vert u\Vert_{H^{1}(\Omega)}^{2}%
\]
Because we assumed that the boundary $\partial\Omega$ is at least $C^{2}$,
regularity theory shows that a solution $u\in H^{1}(\Omega)$ of the
variational problem
\begin{equation}
\mathcal{A}(u,v)=(f,v)\quad\text{for all }v\in H^{1}(\Omega)\label{eq1003}%
\end{equation}
fulfills $u\in D_{L}$; see again \cite{LadUra} or \cite{Triebel}. So the
problems (\ref{eq1002}) and (\ref{eq1003}) are equivalent.

Instead of (\ref{eq1001}) we consider the equivalent variational problem to
find $u\in H^{1}(\Omega)$ which solves
\[
\mathcal{A}(u,v)=(\lambda+1)\int_{\Omega}u(s)v(s)\,ds\quad\text{for all }v\in
H^{1}(\Omega)
\]
and this is equivalent to
\begin{equation}
\int_{\Omega}\nabla_{s}u(s)\nabla_{s}v(s)\,ds=\lambda\,\int_{\Omega
}u(s)v(s)\,ds\quad\text{for all }v\in H^{1}(\Omega)\label{eq1004}%
\end{equation}
Equation (\ref{eq1004}) is the starting point for our numerical approximation
scheme, see also \ref{e24}. First we transfer equation (\ref{eq1004}) to an
equation on the domain $B_{3}$ with the help of a transformation $\Phi
:B_{3}\mapsto\Omega$. So (\ref{eq1004}) becomes
\begin{equation}%
\begin{array}
[c]{l}%
{\displaystyle\int_{B_{3}}}
\nabla_{x}\widetilde{u}(x)\widetilde{A}(x)\nabla_{x}\widetilde{v}%
(x)|\det(J(x)|\,dx\smallskip\\
\left.  \quad\right.  =\lambda%
{\displaystyle\int_{B_{3}}}
\widetilde{u}(x)\widetilde{v}(x)|\det(J(x)|\,dx\quad\text{for all
}\widetilde{v}\in H^{1}(B_{3})
\end{array}
\label{eq1005}%
\end{equation}
where $\widetilde{A}(x)=J(x)^{-1}J(x)^{-T}$; see (\ref{e3a})--(\ref{e3c}) for
the definition of the functions and $J(x)$. According to Section \ref{sec2.2}
we need a sequence of subspaces $\widetilde{\mathcal{X}}_{n}\subset
\widetilde{\mathcal{X}}_{n+1}\subset H^{1}(B_{3})$ with
\[
\overline{\bigcup_{n=1}^{\infty}\widetilde{\mathcal{X}}_{n}}=H_{1}(B_{3})
\]
Because there are no boundary conditions imposed on $H^{1}(B_{3})$ we can use
\[
\mathcal{X}_{n}=\{p(x)\;|\;p\in\Pi_{n}\}
\]
where $\Pi_{n}$ is the space of polynomials in $3$ variables of degree $n$ or
less. As a basis we choose
\[
\{\widetilde{\varphi}_{i}(x)\;|\;i=1,\ldots,N_{n}\},\quad N_{n}={\binom
{n+3}{3}}%
\]
where $\left\{  \widetilde{\varphi}_{i}\right\}  $ is an enumeration of the
orthogonal basis $\left\{  \widetilde{\varphi}_{m,j,\beta}\right\}  $ given in
(\ref{e1001}). To approximate the solutions $\widetilde{u}(x)$ of
(\ref{eq1005}) we use
\[
\widetilde{u}_{n}^{(i)}(x)=\sum_{j=1}^{N_{n}}\alpha_{j}^{(i)}%
\widetilde{\varphi}_{j}(x),:i=,\ldots,N_{n}%
\]
and the coefficients $\alpha_{j}^{(i)}$ for the eigenvalue approximation
$\lambda_{n}^{(i)}$ are given as solutions of the finite eigenvalue problem%
\begin{equation}%
\begin{array}
[c]{l}%
{\displaystyle\sum\limits_{j=1}^{N_{n}}}
\left(
{\displaystyle\int_{B_{3}}}
\nabla_{x}\widetilde{\varphi}_{j}(x)\widetilde{A}(x)\nabla_{x}%
\widetilde{\varphi}_{k}(x)|\det(J(x))|\,dx\right)  \alpha_{j}^{(i)}%
\smallskip\\
\left.  \quad\right.  =\lambda_{n}^{(i)}%
{\displaystyle\sum\limits_{j=1}^{N_{n}}}
\left(
{\displaystyle\int_{B_{3}}}
\widetilde{\varphi}_{j}(x)\widetilde{\varphi}_{k}(x)|\det(J(x))|\,dx\right)
\alpha_{j}^{(i)},\quad k=1,\dots,N_{n}%
\end{array}
\label{eq1006}%
\end{equation}
The functions $\nabla_{x}\widetilde{\varphi}_{j}(x)$ can be calculated
explicitly and all integrals in formula (\ref{eq1006}) are approximated by the
quadrature formula (\ref{e1004}) with $q=n$. The convergence analysis of
Section \ref{sec2.3} can be used without any modifications.%

\begin{figure}[ptb]%
\centering
\includegraphics[
height=3in,
width=3.9998in
]%
{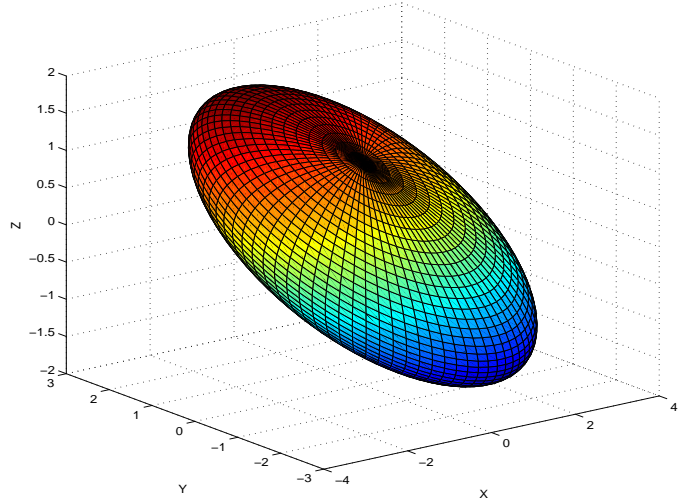}%
\caption{The boundary of $\Omega_{1}$}%
\label{figure_ellipsoid}%
\end{figure}
\begin{figure}[tb]%
\centering
\includegraphics[
height=3in,
width=3.9998in
]%
{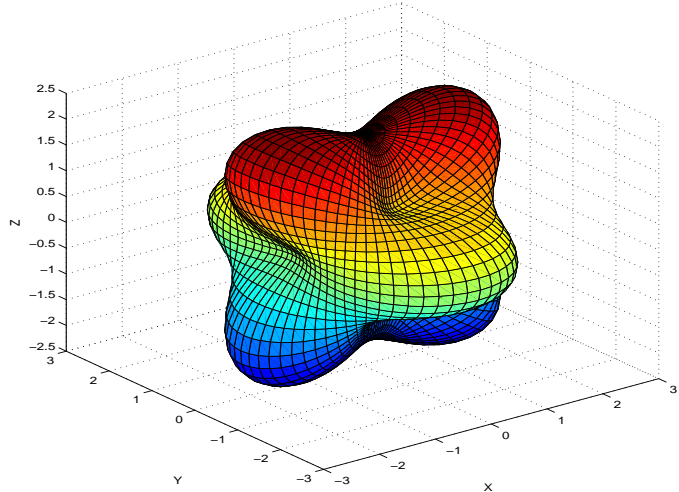}%
\caption{A view of $\partial\Omega_{2}$}%
\label{shape_omega2}%
\end{figure}

To test our method we use two different domains. Let $B_{3}$ denote the closed
unit ball in $\mathbb{R}^{3}$. The domain $\Omega_{1}=\Phi_{1}(B_{3})$ is
given by
\[
s=\Phi_{1}(x)\equiv\left(
\begin{array}
[c]{c}%
x_{1}-3x_{2}\\
2x_{1}+x_{2}\\
x_{1}+x_{2}+x_{3}%
\end{array}
\right)
\]
so $B_{3}$ is transformed to an ellipsoid $\Omega_{1}$; see Figure
\ref{figure_ellipsoid}. The domain $\Omega_{2}$ is given by
\begin{equation}
\Phi_{2}\left(
\begin{array}
[c]{r}%
\rho\\
\phi\\
\theta
\end{array}
\right)  =\left(
\begin{array}
[c]{c}%
(1-t(\rho))\rho+t(\rho)S(\phi,\theta)\\
\phi\\
\theta
\end{array}
\right) \label{eq1007}%
\end{equation}
where we used polar coordinates $(\rho,\phi,\theta)\in\lbrack0,1]\times
\lbrack0,2\pi]\times\lbrack0,\pi]$ to define the mapping $\Phi_{2}$. Here the
function $S:S^{2}=\partial B_{3}\mapsto(1,\infty)$ is a function which
determines the boundary of a star shaped domain $\Omega_{2}$. The restriction
$S(\phi,\theta)>1$ guarantees that $\Phi_{2}$ is injective, and this can
always be assumed after a suitable scaling of $\Omega_{2}$. For our numerical
example we use
\[
S(\theta,\phi)=2+\frac{3}{4}\cos(2\phi)\sin(\theta)^{2}(7\cos(\theta)^{2}-1)
\]
Finally the function $t$ is defined by
\[
t(\rho)\equiv\left\{
\begin{array}
[c]{cc}%
0, & 0\leq\rho\leq\frac{1}{2},\\
2^{5}(\rho-\frac{1}{2})^{5}, & \frac{1}{2}<\rho\leq1.
\end{array}
\right.
\]
where the exponent $5$ implies $\Phi_{2}\in C^{4}(B_{1}(0))$. See
\cite{ach2009} for a more detailed description of $\Phi_{2}$; one perspective
of the surface $\Omega_{2}$ is shown in Figure \ref{shape_omega2}.%

\begin{table}[tbp] \centering
\label{table1000}\caption{Numerical results for $\Omega_1$, $h=0.0001$ to
approximate $R^{(i)}_n$}%
\begin{tabular}
[c]{|c|c|c|c|c|c|c|c|}\hline
$n$ & $N_{n}$ & $|\lambda_{n}^{(1)}-\lambda_{15}^{(1)}|$ & $|\lambda_{n}%
^{(2)}-\lambda_{15}^{(2)}|$ & $\angle(u_{n}^{(1)},u_{15}^{(1)})$ &
$\angle(u_{n}^{(2)},u_{15}^{(2)})$ & $R_{n}^{(1)}$ & $R_{n}^{(2)}$\\\hline
$1$ & $4$ & $4.26E-2$ & $1.00E-1$ & $9.93E-2$ & $1.31E-1$ & $1.45E-1$ &
$2.63E-1$\\\hline
$2$ & $10$ & $4.26E-2$ & $1.00E-1$ & $9.93E-2$ & $1.31E-1$ & $1.45E-1$ &
$2.63E-1$\\\hline
$3$ & $20$ & $1.42E-4$ & $5.67E-4$ & $5.47E-3$ & $1.01E-2$ & $2.28E-2$ &
$5.15E-2$\\\hline
$4$ & $35$ & $1.42E-4$ & $5.67E-4$ & $5.47E-3$ & $1.01E-2$ & $2.28E-2$ &
$5.15E-2$\\\hline
$5$ & $56$ & $1.06E-7$ & $8.38E-7$ & $1.04E-4$ & $2.72E-4$ & $1.22E-3$ &
$3.54E-3$\\\hline
$6$ & $84$ & $1.06E-7$ & $8.38E-7$ & $1.04E-4$ & $2.72E-4$ & $1.22E-3$ &
$3.54E-3$\\\hline
$7$ & $120$ & $2.53E-11$ & $4.31E-10$ & $1.24E-6$ & $4.85E-6$ & $3.02E-5$ &
$1.08E-4$\\\hline
$8$ & $165$ & $2.53E-11$ & $4.31E-10$ & $1.24E-6$ & $4.85E-6$ & $3.02E-5$ &
$1.08E-4$\\\hline
$9$ & $220$ & $2.22E-11$ & $6.78E-14$ & $0$ & $6.32E-8$ & $4.25E-7$ &
$1.80E-6$\\\hline
$10$ & $286$ & $4.47E-11$ & $1.81E-13$ & $0$ & $5.77E-8$ & $4.21E-7$ &
$1.80E-6$\\\hline
$11$ & $364$ & $1.84E-13$ & $5.19E-13$ & $0$ & $0$ & $1.48E-8$ &
$4.15E-8$\\\hline
$12$ & $455$ & $2.07E-13$ & $1.18E-13$ & $0$ & $0$ & $2.21E-9$ &
$1.88E-8$\\\hline
$13$ & $560$ & $1.52E-13$ & $1.91E-13$ & $0$ & $0$ & $5.81E-9$ &
$3.43E-8$\\\hline
$14$ & $680$ & $4.64E-13$ & $5.56E-14$ & $0$ & $0$ & $1.21E-8$ &
$4.26E-8$\\\hline
\end{tabular}%
\end{table}%
%

\begin{table}[tbp] \centering
\label{table1001}\caption{Numerical results for $\Omega_2$}%
\begin{tabular}
[c]{|c|c|c|c|c|c|}\hline
$n$ & $N_{n}$ & $|\lambda_{n}^{(1)}-\lambda_{15}^{(1)}|$ & $|\lambda_{n}%
^{(2)}-\lambda_{15}^{(2)}|$ & $\angle(u_{n}^{(1)},u_{15}^{(1)})$ &
$\angle(u_{n}^{(2)},u_{15}^{(2)})$\\\hline
$1$ & $4$ & $3.60E-1$ & $3.21E-1$ & $2.86E-1$ & $3.31E-1$\\\hline
$2$ & $10$ & $3.60E-1$ & $3.21E-1$ & $2.86E-1$ & $3.31E-1$\\\hline
$3$ & $20$ & $8.16E-2$ & $8.53E-2$ & $8.32E-2$ & $1.05E-1$\\\hline
$4$ & $35$ & $8.16E-2$ & $8.53E-2$ & $8.32E-2$ & $1.05E-1$\\\hline
$5$ & $56$ & $1.99E-2$ & $2.27E-2$ & $2.84E-2$ & $3.11E-2$\\\hline
$6$ & $84$ & $1.99E-2$ & $2.27E-2$ & $2.84E-2$ & $3.11E-2$\\\hline
$7$ & $120$ & $1.48E-2$ & $1.56E-2$ & $2.49E-2$ & $2.71E-2$\\\hline
$8$ & $165$ & $1.48E-2$ & $1.56E-2$ & $2.49E-2$ & $2.71E-2$\\\hline
$9$ & $220$ & $4.77E-3$ & $5.96E-3$ & $1.14E-2$ & $1.47E-2$\\\hline
$10$ & $286$ & $4.77E-3$ & $5.96E-3$ & $1.14E-2$ & $1.47E-2$\\\hline
$11$ & $364$ & $8.34E-4$ & $1.28E-3$ & $3.25E-3$ & $4.76E-3$\\\hline
$12$ & $455$ & $8.34E-4$ & $1.28E-3$ & $3.25E-3$ & $4.76E-3$\\\hline
$13$ & $560$ & $1.88E-4$ & $2.55E-4$ & $1.33E-3$ & $1.62E-3$\\\hline
$14$ & $680$ & $1.88E-4$ & $2.55E-4$ & $1.33E-3$ & $1.62E-3$\\\hline
\end{tabular}%
\end{table}%
%

\begin{figure}[tb]%
\centering
\includegraphics[
height=3in,
width=3.9998in
]%
{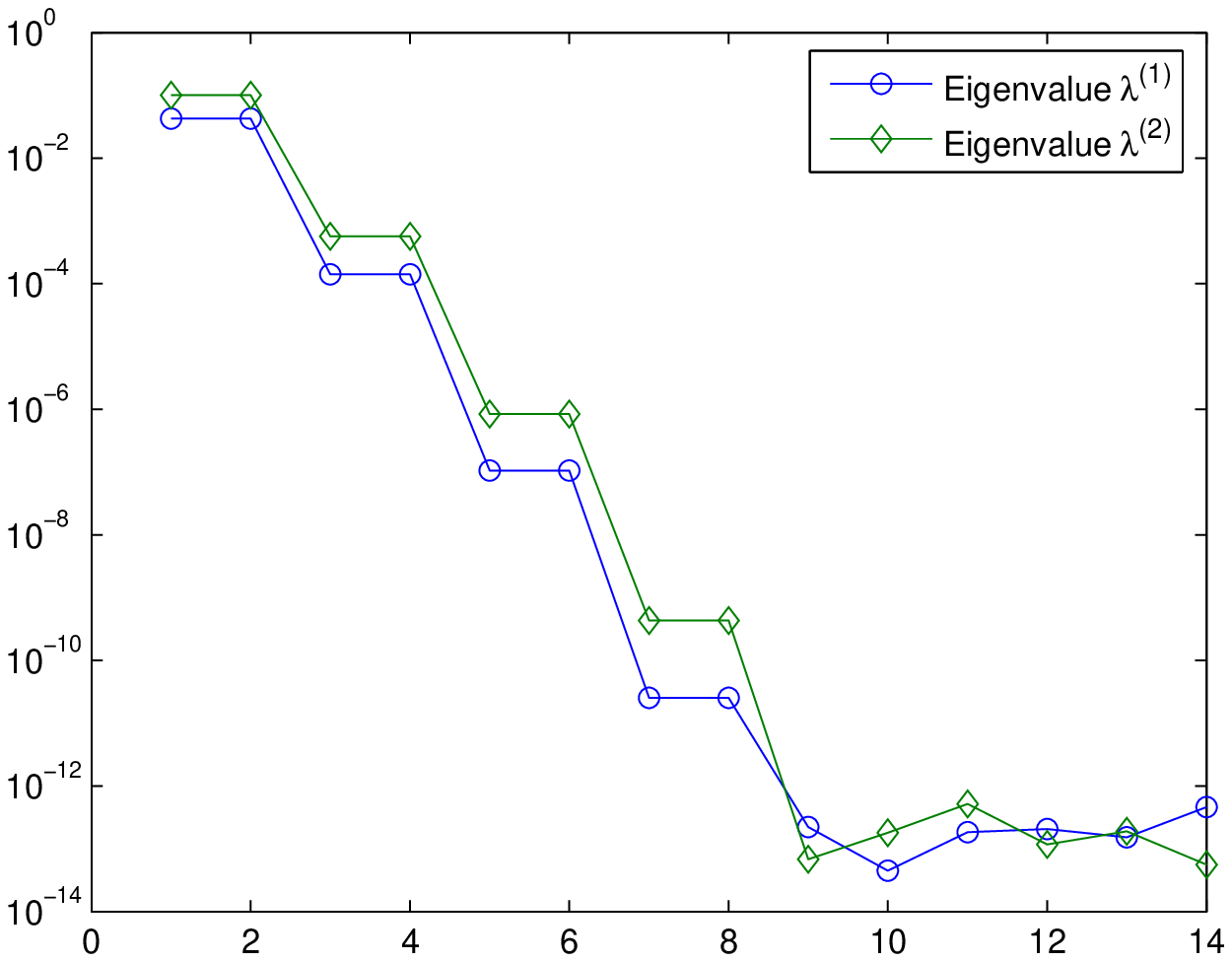}%
\caption{$\Omega_{1}$: errors $\left\vert \lambda_{15}^{(i)}-\lambda_{n}%
^{(i)}\right\vert $ for the calculation of the first two eigenvalues
$\lambda^{(i)}$}%
\label{figure_omega3lambda}%
\end{figure}
\begin{figure}[ptb]%
\centering
\includegraphics[
height=3in,
width=3.9998in
]%
{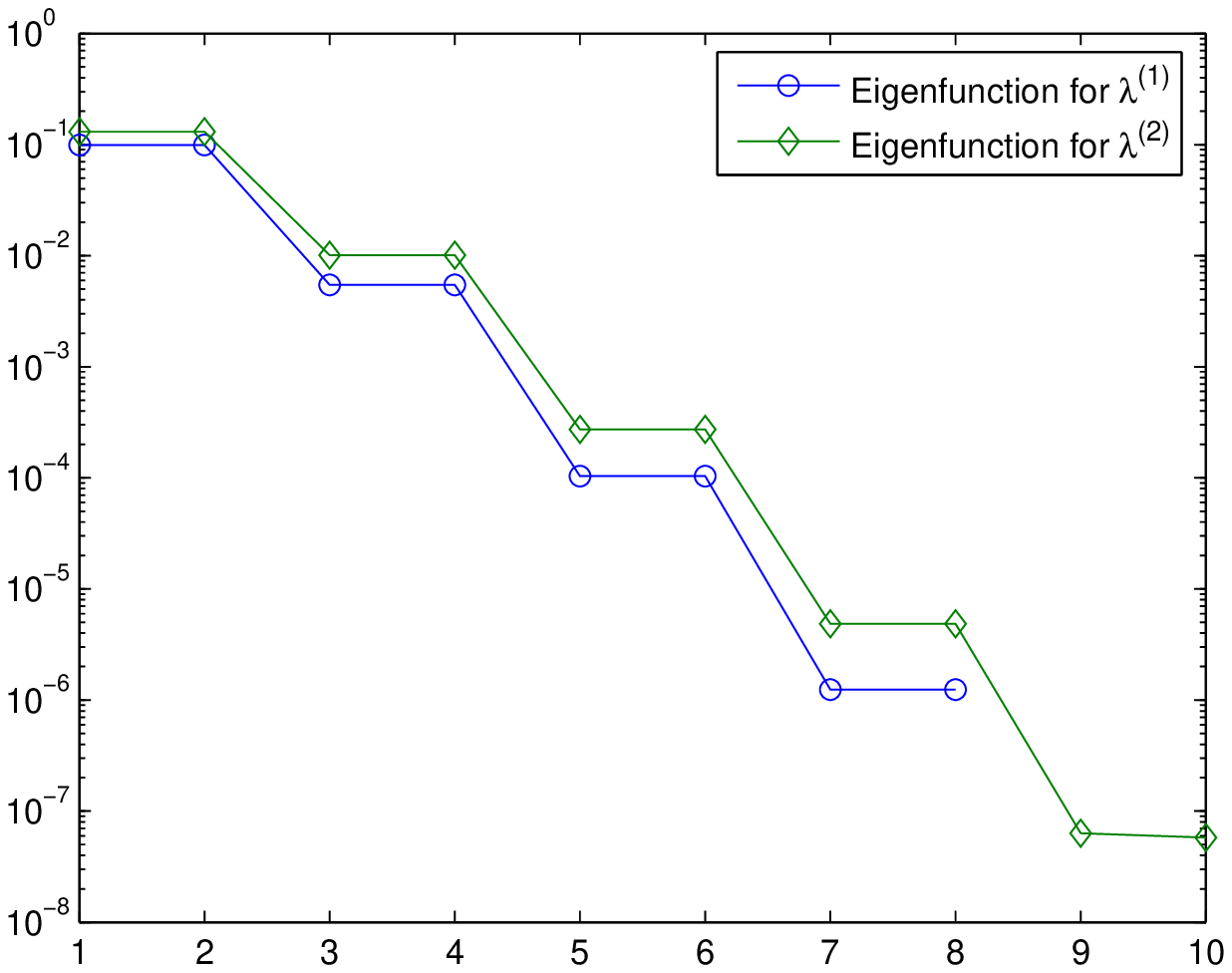}%
\caption{$\Omega_{1}$: angles $\angle\left(  u_{n}^{(i)},u_{15}^{(i)}\right)
$ between the approximate eigenfunction $u_{n}^{(i)}$ and our most accurate
approximation $u_{15}^{(i)}\approx u^{(i)}$.}%
\label{figure_omega3angle}%
\end{figure}
\begin{figure}[tb]%
\centering
\includegraphics[
height=3in,
width=3.9998in
]%
{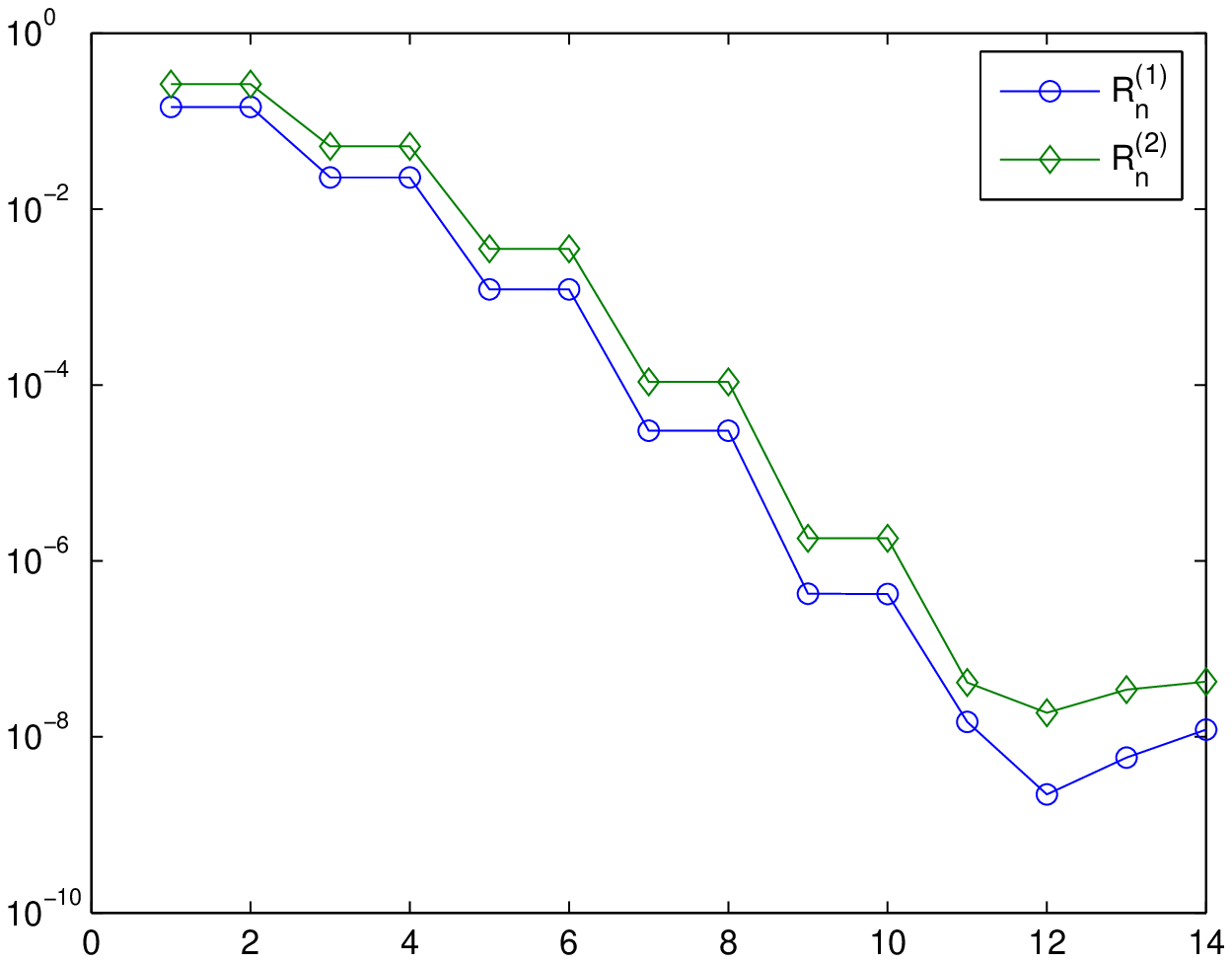}%
\caption{$\Omega_{1}$: errors $\left\vert -\Delta u_{n}^{(i)}\left(  s\right)
-\lambda_{n}^{(i)}u_{n}^{(i)}\left(  s\right)  \right\vert $}%
\label{figure_omega3delta}%
\end{figure}

For each domain we calculate the approximate eigenvalues $\lambda_{n}^{(i)}$,
$\lambda_{n}^{(0)}=0<\lambda_{n}^{(1)}\leq\lambda_{n}^{(2)}\leq\ldots$ and
eigenfunctions $u_{n}^{(i)}$, $i=1,\ldots,N_{n}$, for the degrees
$n=1,\ldots,15$ (here we do not indicate dependence on the domain $\Omega$).
To analyze the convergence we calculate several numbers. First we estimate the
speed of convergence for the first two eigenvalues by calculating
$|\lambda_{15}^{(i)}-\lambda_{n}^{(i)}|$, $i=1,2$, $n=1,\ldots,14$. Then to
estimate the speed of convergence of the eigenfunctions we calculate the angle
(in $L^{2}(\Omega)$) between the current approximation and the most accurate
approximation $\angle(u_{n}^{(i)},u_{15}^{(i)})$, $i=1,2$, $n=1,\ldots,14$.
Finally, an independent estimate of the quality of our approximation is given
by
\[
R_{n}^{(i)}\equiv|-\Delta u_{n}^{(i)}(s)-\lambda_{n}^{(i)}u(s)|
\]
where we use only one $s\in\Omega$, given by $\Phi(1/10,1/10,1/10)$. To
approximate the Laplace operator we use a second order difference scheme with
$h=0.0001$ for $\Omega_{1}$ and $h=0.01$ for $\Omega_{2}$. The reason for the
latter choice of $h$ is that our approximations for the eigenfunctions on
$\Omega_{2}$ are only accurate up three to four digits, so if we divide by
$h^{2}$ the discretization errors are magnified to the order of $1$.

The numerical results for $\Omega_{1}$ are given in table \ref{table1000}. The
graphs in Figures \ref{figure_omega3lambda}--\ref{figure_omega3delta}, seem to
indicate exponential convergence. For the graphs of $\angle(u_{n}^{(i)}%
,u_{15}^{(i)})$, see Figure \ref{figure_omega3angle}. We remark that we use
the function $\arccos(x)$ to calculate the angle, and for $n\approx9$ the
numerical calculations give $x=1$, so the calculated angle becomes $0$. For
the approximation of $R_{n}^{(i)}$ one has to remember that we use a
difference method of order $O(h^{2})$ to approximate the Laplace operator, so
we can not expect any result better than $10^{-8}$ if we use $h=0.0001$.

\begin{figure}[tb]%
\centering
\includegraphics[
height=3in,
width=3.9998in
]%
{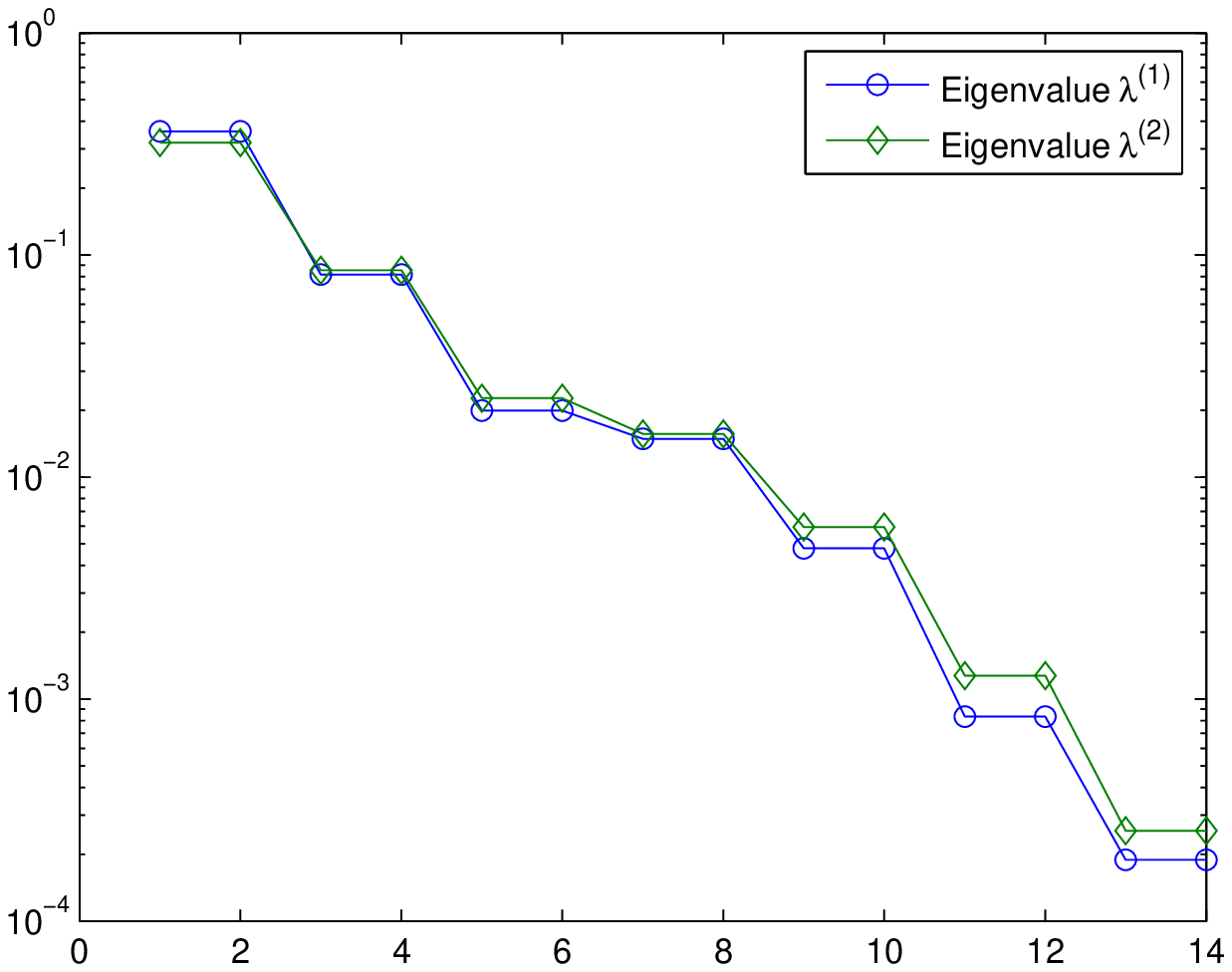}%
\caption{$\Omega_{2}$: errors $\left\vert \lambda_{15}^{(i)}-\lambda_{n}%
^{(i)}\right\vert $ for the calculation of the first two eigenvalues
$\lambda^{(i)}$}%
\label{figure_omega4lambda}%
\end{figure}

\begin{figure}[tb]%
\centering
\includegraphics[
height=3in,
width=3.9998in
]%
{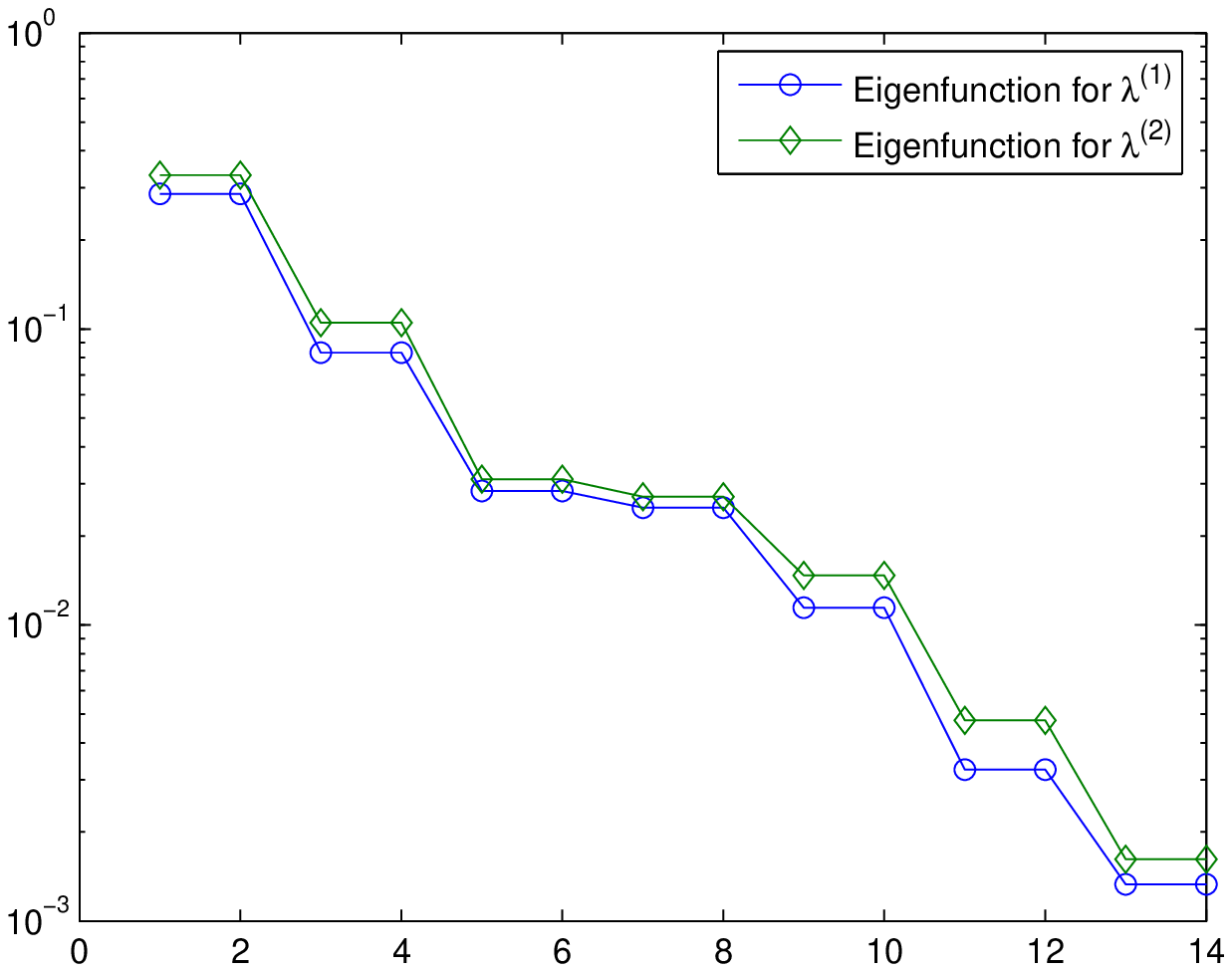}%
\caption{$\Omega_{2}$: angles $\angle\left(  u_{n}^{(i)},u_{15}^{(i)}\right)
$ between the approximate eigenfunction $u_{n}^{(i)}$ and our most accurate
approximation $u_{15}^{(i)}\approx u^{(i)}$.}%
\label{figure_omega4angle}%
\end{figure}

As we expect, the approximations for $\Omega_{2}$ with the transformation
$\Phi_{2}$ present a bigger problem for our method. Still from the graphs in
Figure \ref{figure_omega4lambda} and \ref{figure_omega4angle} we might infer
that the convergence is exponential, but with a smaller exponent than for
$\Omega_{1}$. Because $\Phi_{2}\in C^{4}(B_{3})$ we know that the transformed
eigenfunctions on $B_{3}$ are in general only $C^{4}$, so we can only expect a
convergence of $O(n^{-4})$. The values of $n$ which we use are too small to
show what we believe is the true behavior of the $R_{n}^{(i)}$, although the
values for $n=10\ldots14$ seem to indicate some convergence of the type we
would expect.

The poorer convergence for $\Omega_{2}$ as compared to $\Omega_{1}$
illustrates a general problem. When defining a surface $\partial\Omega$ by
giving it as the image of a 1-1 mapping from the unit sphere $S^{2}$ into
$\mathbb{R}^{3}$, how does one extend it to a smooth mapping from the unit
ball to $\Omega$? The mapping in (\ref{eq1007}) is smooth, but it has large
changes in its derivatives, and this affects the rate of convergence of our
spectral method. We are working at present on this problem, developing a
numerical method to find a well-behaved polynomial mapping $\Phi$ when given
only its restriction to $S^{2}$.%

\end{document}